\documentclass[12pt]{amsart}

\usepackage{amssymb,latexsym,amsmath}
\usepackage{verbatim,fullpage}

\input {cyracc.def}
\font\ququ=cmr10 scaled \magstep1
\font\tencyr=wncyr8 scaled \magstep1

\def\rus{\tencyr\cyracc}

\newcommand{\re}[1]{\textrm  (\ref{#1})}

\renewenvironment{proof}
{\noindent {\sl Proof.}\quad }{\hfill $\square$
\vskip1.1ex\noindent }

\newenvironment{proof*}
{\noindent {\sl Proof.}\quad }{\hfill $\square$}

\renewcommand{\theequation}{\thesection.\arabic{equation}}
\renewcommand{\thesubsubsection}{\theequation.\arabic{subsubsection}}

\catcode`\@=\active
\catcode`\@=11
\def\@eqnnum{\hbox to
.01pt{}\rlap{\hskip-\displaywidth(\mathbf{\theequation})}}

\catcode`\@=12

\newenvironment{s}[1]
{ \vskip1.2ex \refstepcounter{equation}
\noindent {\bf \theequation\quad #1.} \begin{sl}}{\end{sl}
\vskip1.1ex\noindent }

\newenvironment{rem}[1]
{ \vskip1.2ex \refstepcounter{equation}
\noindent {\bf \theequation\enspace {#1}.} }{ \vskip1.1ex\noindent }

\newcommand {\g}{{\frak g}}
\newcommand {\h}{{\frak h}}
\newcommand {\ka}{{\frak k}}
\newcommand {\el}{{\frak l}}

\newcommand {\n}{{\frak n}}

\newcommand {\p}{{\frak p}}
\newcommand {\q}{{\frak q}}

\newcommand {\te}{{\frak t}}
\newcommand {\ut}{{\frak u}}
\newcommand {\z}{{\frak z}}
\newcommand {\sln}{{\frak sl}_{n+1}}
\newcommand {\slv}{{\frak sl}(V)}
\newcommand {\glv}{{\frak gl}(V)}
\newcommand {\spv}{{\frak sp}(V)}
\newcommand {\spn}{{\frak sp}_{2n}}

\newcommand {\esi}{\varepsilon}
\newcommand {\ap}{\alpha}

\newcommand {\lb}{\lambda}

\newcommand {\vp}{\varphi}

\newcommand {\ca}{{\mathcal A}}

\newcommand {\N}{{\mathcal N}}

\newcommand {\co}{{\mathcal O}}
\newcommand {\cs}{{\mathcal S}}

\newcommand {\VV}{{V}}

\newcommand {\md}{/\!\!/}

\newcommand {\isom}{\stackrel{\sim}{\longrightarrow}}
\newcommand {\ad}{{\mathrm{ad\,}}}
\newcommand {\Ad}{{\mathrm{Ad\,}}}

\newcommand {\Aut}{{\mathrm{Aut\,}}}

\newcommand {\hot}{{\mathrm{ht}}}
\newcommand {\ind}{{\mathrm{ind\,}}}

\newcommand {\Lie}{{\mathrm{Lie\,}}}
\newcommand {\Ker}{{\mathrm{Ker\,}}}
\newcommand {\Ima}{{\mathrm{Im\,}}}

\newcommand {\rk}{{\mathrm{rk\,}}}

\newcommand {\trdeg}{{\mathrm{trdeg\,}}}

\newcommand {\tri}{{\frak sl}_2}
\newcommand {\GR}[2]{{\mathrm{{ #1}}}_{#2}}

\newcommand {\ov}{\overline}
\newcommand {\un}{\underline}

\newcommand {\gthe}{\g\langle 1\rangle_{(\theta)}}
\newcommand {\gedva}{\g\langle {\ge }2\rangle}
\newcommand {\vno}[1]{\vskip#1 ex\noindent}
\newcommand {\rar}{\rightarrow}

\newcommand {\beq}{\begin{equation}}
\newcommand {\eeq}{\end{equation}}

\newcommand{\zge}{\z_\g(e)}

\newcommand{\vlb}{\VV_\lb}

\newfam\Bbbfam\newfam\eufam\newfam\Eufam%
\font\Bbbfont=msbm10 scaled 1200%
\font\olala=msam10 scaled 1200%
\font\frak=eufm10 scaled 1400%
\font\Bbbsmallfont=msbm8%
\textfont\Bbbfam=\Bbbfont\scriptfont\Bbbfam=\Bbbsmallfont
\font\euzw=eufm10 scaled 1200%
\font\euac=eufm7 scaled 1200%
\font\euacc=eufm7 scaled 1000%
\textfont\eufam=\euzw\scriptfont\eufam=\euac%
\scriptscriptfont\eufam=\euacc%
\font\Euzw=eufm10 scaled \magstep2%
\font\Euac=eufm7 scaled \magstep2%
\textfont\Eufam=\Euzw\scriptfont\Eufam=\Euac%
\def\frak{\fam\eufam}%
\def\Bbb{\fam\Bbbfam}%

\def\varnothing{\hbox {\Bbbfont\char'077}}
\def\square{\hbox {\olala\char"03}}
\def\bbk{\hbox {\Bbbfont\char'174}}
\def\semidir{\hbox {\olala\char"68}}

\begin{document}
\setlength{\parskip}{2pt plus 4pt minus 0pt}
\hfill {\scriptsize July 15, 2002} \vskip1ex
\vskip1ex

\title[]{
Some amazing properties of spherical nilpotent orbits
}
\author[]
{\sc Dmitri I. Panyushev} 
\maketitle
\begin{center}
{\footnotesize
{\it Independent University of Moscow,
Bol'shoi Vlasevskii per. 11 \\
121002 Moscow, \quad Russia \\ e-mail}: {\tt panyush@mccme.ru }\\
}
\end{center}

\medskip
\smallskip

\section*{Introduction}
\vno{1}%
Let $G$ be a simple algebraic group defined
over an algebraically closed field $\bbk$ of characteristic zero.
Write $\g$ for its Lie algebra.
\\[.6ex]
Let $x\in\g$ be a nilpotent element and $G{\cdot}x\subset \g$
the corresponding nilpotent orbit. The maximal number $m$ such that
$(\ad x)^m\ne 0$ is called the {\it height\/} of
$x$ or of $G{\cdot}x$, denoted
$\hot(x)$. Recall that an irreducible $G$-variety $X$
is called $G$-{\it spherical\/} if a Borel subgroup of $G$ has an open orbit
in $X$. It was shown in \cite{p94} that
$G{\cdot}x$ is $G$-spherical if and only if $(\ad x)^4=0$.
This means that the spherical nilpotent orbits are precisely
the orbits of height 2 and 3.
Unfortunately, whereas my proof for the orbits of height 2 and height $\ge 4$
was completely general, the argument for the orbits of height 3
explicitly used their classification.
In this paper, we give a proof of sphericity that does not rely on the
classification of nilpotent orbits, see Theorem~\ref{height3}.
We begin with some properties of invariants of symplectic
representations. For instance, we prove that {\sf (1)} if
$H\subset Sp(V)$ is an irreducible representation without non-constant
invariants, then $H=Sp(V)$, and {\sf (2)} if $H$ has non-constant invariants,
then it has an invariant of degree 4.
Applying these results to nilpotent orbits, we prove that
the centraliser $\z_\g(x)$ has a rather specific
structure whenever  $\hot(x)$ is odd. From this description,
we deduce a conceptual proof of sphericity in case $\hot(x)=3$.
As another application we compute the index of $\z_\g(x)$.
It will be shown that $\ind\z_\g(x)=\rk\g$, if $\hot(x)=3$.
This confirms Elashvili's conjecture for such $x$ (see \cite[Sect.\,3]{p02}
about this conjecture). In Section~\ref{special},
we prove that if $\theta$, the highest root of
$\g$, is fundamental, then $\g$ always has a specific orbit of height 3,
which is denoted by $\Bbb O$.
This orbit satisfies several arithmetical relations. Namely,
if $\g=
\bigoplus_{-3\le i\le 3}\g\langle i\rangle $ is the
$\Bbb Z$-grading associated with an $\tri$-triple containing $e\in\Bbb O$, then
$\dim\g\langle 3\rangle=2$ and $\dim\g\langle 1\rangle =
2\,\dim\g\langle 2\rangle $. Furthermore, the weighted Dynkin diagram of
$\Bbb O$ can explicitly be described.
Let $\beta$ be the unique simple root that is not orthogonal
to $\theta$ and let $\{\ap_i\}$ be all simple roots adjacent to $\beta$ on
the Dynkin diagram of $\g$. Then one has to put `1' at all $\ap_i$'s and
`0' at all other simple roots (Theorem~\ref{char_3}).
It is curious that these properties of $\Bbb O$ enables us to give
an intrinsic construction of $\GR{G}{2}$-grading in each simple $\g$
whose highest root is fundamental.
\\[.6ex]
In Section 5, the problem of computing
the algebra of covariants on a nilpotent orbit is discussed. Let
$\g=\underset{i}{\bigoplus}
\g\langle i\rangle$ be the $\Bbb Z$-grading associated with $e$. Set
$L=\exp\g\langle 0\rangle\subset G$. Let $\ut_-\oplus\te\oplus\ut_+$
be a triangular decomposition of $\g$ such that $\te\subset\g\langle 0\rangle$
and $\ut_+\supset \g\langle{\ge}1\rangle$. Put
$U_-=\exp(\ut_-)$ and $U(L)_-=\exp(\ut_-\cap\g\langle 0\rangle)$. These are
maximal unipotent subgroups in $G$ and $L$, respectively. We show that
there is an injective homomorphism
$   \tau^o : \bbk[G{\cdot}e]^{U_-} \to \bbk[\gedva]^{U(L)_-}$,
which is birational, i.e., induces an isomorphism of the quotient fields
(Theorem~\ref{restrict}).
This result can be restated in the ``dominant'' form, when the
unipotent groups in question are replaced by the opposite ones and the algebra
of functions on $\gedva$ is replaced by the symmetric algebra. Namely,
the natural homomorphism
$\hat\tau^o : \bbk[G{\cdot}e]^{U} \to \cs(\gedva)^{U(L)}$
is injective and birational.
A complete answer is obtained if $\hot(e)=2$.
In this case, it is proved that $\tau^o$ is an isomorphism and
$\bbk[\ov{G{\cdot}e}]$ is a free $\bbk[\ov{G{\cdot}e}]^{U}$-module,
see Theorem~\ref{height2}. It is also shown how one can
quickly determine the monoid of highest weights of all simple $G$-modules in
the algebra of regular functions $\bbk[\ov{G{\cdot}e}]$.
The case of height~3 is more complicated. However,
we have a general conjecture describing the range of $\tau^o$, if $G{\cdot}e$
is spherical, see \re{conj1}.
If true, this conjecture allows us to obtain a complete
description of $\bbk[G{\cdot}e]^U$.
The agreement of this description with known results for orbits in
exceptional Lie algebras
is in my view a significant
evidence for the validity of Conjecture~\ref{conj1}.
One of the advantages of this method is that no computer calculations is
needed. Furthermore, our
approach coupled with results of \cite{e8} yields a description of
the algebra of polynomial functions on the model orbit in
$\GR{E}{8}$ as {\it graded\/} $G$-module.
Our computations for the spherical nilpotent orbits
are gathered in two Tables in Section~\ref{tables}.

{\sl Notation.} Algebraic groups are denoted by capital Latin characters
and their Lie algebras are usually denoted by the
corresponding small Gothic characters; $H^0$ is the identity component of
an algebraic group $H$. If $H\hookrightarrow GL(V)$ and $v\in V$,
then $H{\cdot}v$ is the orbit and $H_v$ is the isotropy group of $v$;
$\h_v$ is the stationary subalgebra of $v$ in $\h=\Lie H$.

{\small {\bf Acknowledgements.}
This research was supported in part by
R.F.B.I. Grant 01--01--00756.
}

\section{On invariants of symplectic representations}
\label{inv_symp}
\noindent
Let $H$ be a reductive group with a maximal torus $T_H$.
Fix also a triangular decomposition of the Lie algebra:
$\h=\n_+\oplus\te_H\oplus \n_-$. Write $\vlb$ for the simple $H$-module
with highest weight $\lb$. We begin with recalling a standard fact on
finite-dimensional representations.

\begin{s}{Lemma}   \label{vesa}
Let $\VV_{\lb_1}$ and $\VV_{\lb_2}$ be two simple $H$-modules.
If $\mu_i$ is an arbitrary weight of $\VV_{\lb_i}$ with respect to $T_H$,
then $\mu_1+\mu_2$ is a weight of
$\VV_{\lb_1+\lb_2}$.
\end{s}%
We are going to prove a general result concerning symplectic
representations without invariants.
The reader may observe that the next proof
can be made shorter, if one invokes the classification
of simple algebraic groups. Since our intention is to argue in a
classification-free way, we prefer to give a longer  but ``pure"
argument.

\begin{s}{Theorem} \label{main1}
Let $H\to Sp(\VV)$ be a faithful symplectic representation of a
reductive group $H$.
Suppose $\bbk[\VV]^H=\bbk$. Then there exists a decomposition
$\VV=\oplus_{i=1}^l\VV_i$ such that
$H=\prod_{i=1}^l Sp(\VV_i)$.
\end{s}\begin{proof*}
$1^o$. Without loss of generality, we may assume that $H$ is connected.
Indeed, if the conclusion holds for $H^0$, then this automatically
implies that $H=H^0$.
Since $(\cs^2\VV)^H=0$ and $\cs^2\VV\simeq\spv$, we see that the centre
of $H$ is finite. Hence $H$ is semisimple.
Let $<\ ,\ >$ denote a skew-symmetric $H$-invariant bilinear form on
$\VV$.

$2^o$. Assume first that $\VV=\vlb$ is a {\it simple\/} $H$-module.
Let $\Delta$ be the root system with respect to $T_H$ and
$\Delta^+$ the subset of positive roots corresponding to $\n_+$.
Let $v_\lb$ (resp. $v_{-\lb}$) be a highest (resp. lowest) weight vector
in $\vlb$.  Let $L_\lb$ denote the standard Levi subgroup in
the parabolic subgroup that stabilises the line $\bbk v_\lb$.
Set $L'_\lb=L_\lb\cap H_{v_\lb}$. It is a codimension-1 reductive
subgroup of $L_\lb$.
Consider the $H$-orbit of $w=v_\lb+v_{-\lb}$.
Because $\bbk[\VV]^H=\bbk$, each $H$-orbit is unstable, i.e.,
its closure contains the origin. In particular,
$\ov{H{\cdot}w}\ni 0$. It is clear that
$L'_\lb\subset H_w$. By Luna's criterion \cite[Cor.\,2]{luna}, we have
\[
\ov{H{\cdot}w}\ni 0  \Leftrightarrow \ov{N_H(L'_\lb){\cdot}w}\ni 0 .
\]
As $\rk L'_\lb=\rk H-1$, there are only two possibilities for
$N_H(L'_\lb)^0$. It is either $T_H$ or a subgroup of semisimple
rank 1. Obviously, $T_H{\cdot}w$ is closed. Therefore the second possibility
must occur, i.e., $N_H(L'_\lb)$ is locally isomorphic to $SL_2 \times \{
\textrm{torus}\}$. Moreover, the difference $\lb-(-\lb)=2\lb$ must be a
root of $SL_2$ and hence of $H$.
Thus, our intermediate conclusion is: {\it
If\/ $\VV_\lb$ is a simple symplectic $H$-module without non-constant
invariants, then $2\lb=:\ap\in\Delta^+$.}

The previous conclusion shows that only one simple component of
$H$ can act non-trivially on $\vlb$. Hence $H$ is a {\it simple\/}
algebraic group.
Since $\ap$ is a dominant root of $H$, it is either the highest
(long) root or the highest short root.
Assume $\ap$ is short, and let $\beta$ be the highest root. Then
$(\ap,\beta^\vee)=1$ and hence $(\lb,\beta^\vee)=1/2$, which is impossible.
Hence  $\ap$ is the highest root in $\Delta^+$.
As $\lb=\frac{1}{2}\ap$, we see that $(\lb,\nu^\vee)\le 1$ for any
$\nu\in\Delta^+$, i.e., $\lb$ is minuscule. Therefore, all $T_H$-weights
in $\vlb$ are simple and nonzero.
We have $\h=\VV_{2\lb}$, and
the embedding of Lie algebras $\h\subset{\frak sp}(\vlb)$
is nothing but the embedding $\VV_{2\lb}\subset \cs^2\vlb$.

Let $T$ be a maximal torus of $Sp(\vlb)$ such that
$T_H\subset T$. Let $\{\pm\esi_i\mid i=1,\dots,n\}$ be the weights of
$T$ in $\vlb$.
Set $\mu_i=\esi_i\mid_{T_H}$.
Since $\lb$ is minuscule, we have $\mu_i\ne 0$ and
$\pm\mu_i\pm\mu_j\ne 0$ ($i\ne j$). It follows from Lemma~\ref{vesa} that
$\pm\mu_i\pm\mu_j$ ($i\ne j)$ and $\pm 2\mu_i$ are the nonzero
weights of $\VV_{2\lb}$,
i.e., the roots of $\h$. This means that, for any two different weights
of $\vlb$, their difference is a root of $\h$. It follows that
$\h{\cdot}v_{\lb}=\vlb$, i.e. the $H$-orbit of highest weight vectors
is dense in $\vlb$. The theory of orbits of highest weight vectors developed
in \cite{vp72} says that in this situation $\cs^n\vlb$ is a simple
$H$-module for any $n$. In particular,
$\VV_{2\lb}=\cs^2\vlb$ and hence $\h={\frak sp}(\vlb)$.

$3^o$. Assume now that $\VV$ is a {\it reducible\/} $H$-module, and let
$\VV=\VV_1\oplus\ldots\oplus
\VV_l$ be a decomposition into the sum of irreducibles.
For $v=(v_1,v_2,\dots,v_l)\in\VV$, define the $H$-invariant polynomial
function $f_{ij}$ by the formula $f_{ij}(v)={<}v_i,v_j{>}$. By the
assumption, $f_{ij}\equiv 0$. Hence the restriction of
$<\ ,\ >$ to each $\VV_i$ is non-degenerate and
$H\subset Sp(\VV_1)\times\ldots\times Sp(\VV_l)$.
From the previous part of the proof, it follows that projection of
$H$ to the $i$-th factor equals $Sp(V_i)$.
If $H$ has a simple factor that is embedded diagonally in the product
$Sp(V_i)\times Sp(V_j)$ (when $\dim V_i=\dim V_j$), then
$\bbk[V_i\oplus V_j]^H\ne \bbk$. Therefore $H$ cannot have diagonally
embedded factors and hence $H=Sp(\VV_1)\times\ldots\times Sp(\VV_l)$.
\end{proof*}%
\begin{s}{Corollary {\ququ (of the Proof)}}  \label{w}
Let $\vlb$ be a simple symplectic $H$-module. The following three
conditions are equivalent: a) $2\lb$ is a root of $\h$;
b) the orbit $H{\cdot}(v_\lb+v_{-\lb})$ is unstable;
c)  $H=Sp(\vlb)$.
\end{s}%
This Corollary has a useful complement.

\begin{s}{Proposition}  \label{degree-4}
Suppose $H\to Sp(\vlb)$ is an irreducible symplectic representation.
If $\bbk[\vlb]^H\ne \bbk$, then $H$ has an invariant of degree 4.
\end{s}\begin{proof}
$1^o$. Let $\Phi$ denote a symmetric non-degenerate $H$-invariant bilinear
form on $\h$ and $<\ ,\ >$ a symplectic form on $\vlb$.
Consider the bilinear mapping $\tilde\psi: \vlb\times \vlb\to \h$ which is
defined by
$\Phi(\tilde\psi(v,w),s)=<s{\cdot}v,w>$, $s\in\h$, $v,w\in\vlb$.
We also need the quadratic mapping $\psi: \vlb\to \h$,
$\psi(v):=\tilde\psi(v,v)$. Define the $H$-invariant polynomial of degree 4 by
$F(v)=\Phi(\psi(v),\psi(v))$. The problem is to prove that $F\not\equiv 0$.
Actually, we show that $F(v_\lb+v_{-\lb})\ne 0$.

$2^o$. By Corollary~\ref{w}, we have $2\lb\not\in \Delta$.
Therefore $\h{\cdot}v_\lb\not\ni v_{-\lb}$.
Then $\Phi(\psi(v_\lb),s)=<s{\cdot}v_\lb, v_\lb>=0$ for any $s\in \h$.
Hence $\psi(v_\lb)=0$, and similarly
$\psi(v_{-\lb})=0$. It then follows from bilinearity that
$\psi(v_\lb+v_{-\lb})=2\tilde\psi(v_\lb,v_{-\lb})$.
Next, $\Phi(\tilde\psi(v_\lb,v_{-\lb}),\n_+)=<\n_+{\cdot}v_\lb,
v_{-\lb}>=0$, and likewise for $\n_-$.
Hence $\tilde\psi(v_\lb,v_{-\lb})\in\te_H$.

$3^o$. Without loss of generality assume that $<v_\lb,v_{-\lb}>=1$.
Then $\Phi(\tilde\psi(v_\lb,v_{-\lb}),s)=(d\lb)(s)$ for any $s\in\te_H$.
It follows that $\tilde\psi(v_\lb,v_{-\lb})\ne 0$ and furthermore
$\tilde\psi(v_\lb,v_{-\lb})$ is orthogonal to $\Ker(d\lb)\subset\te_H$
with respect to $\Phi$. Since $\lb$ is an
algebraic character of $T_H$, the restriction of $\Phi$ to
$\Ker(d\lb)$ is non-degenerate, i.e., $\tilde\psi(v_\lb,v_{-\lb})\not\in
\Ker(d\lb)$. Therefore
\[
0\ne d\lb(\tilde\psi(v_\lb,v_{-\lb}))=
\Phi(\tilde\psi(v_\lb,v_{-\lb}),\tilde\psi(v_\lb,v_{-\lb}))=
\frac{1}{4}F(v_\lb+v_{-\lb})  .
\]
\end{proof}%
{\bf Remarks}. 1. Using Luna's criterion \cite[Cor.\,1]{luna},
one can prove that if
$H\ne Sp(\vlb)$, then the orbit $H{\cdot}(v_\lb+v_{-\lb})$ is closed.
\\[.5ex]
2. $\psi: \vlb\to \h\simeq \h^\ast$ is the moment mapping associated
with the symplectic $H$-variety $\vlb$. Similarly, regarding
$\vlb\times \vlb$ as the cotangent bundle of $\vlb$, one sees
that $\tilde\psi$ is the moment mapping for this $H$-variety.

\noindent

\section{Properties of the $\Bbb Z$-grading associated with an $\tri$-triple}
\label{z-graded}
\setcounter{equation}{0}
\noindent
Let $\g$ be a simple Lie algebra with a fixed triangular decomposition
$\g=\ut_-\oplus\te\oplus\ut_+$ and $\Delta$ the corresponding root system.
The roots of $\ut_+$ are positive. Write $\Delta_+$ (resp. $\Pi$) for
the set of positive (resp. simple) roots; $\theta$ is the
highest root in $\Delta^+$. The fundamental weight corresponding to the
$i$-th simple root is denoted by $\vp_i$.
The Killing form on $\g$ is denoted by $\Phi$, and the induced bilinear form
on $\te^*_{\Bbb Q}$ is denoted by $(\ ,\ )$.
For $x\in\g$, let $Z_G(x)$ and $\z_\g(x)$ denote the
centralisers in $G$ and $\g$, respectively. If $M$ is a subset of $\g$,
then $\z_M(x)=\z_\g(x)\cap M$.
Let $\N \subset\g$ be
the nilpotent cone. By the Morozov-Jacobson theorem, each nonzero
element $e\in\N$ can be
included in an ${\frak sl}_2$-triple $\{e,h,f\}$ (i.e.,
$[e,f]=h,\ [h,e]=2e,\ [h,f]=-2f$). The semisimple element $h$, which is called
a {\it characteristic} of $e$,
determines a $\Bbb Z$-grading in $\g$:
\[
 \g=\bigoplus_{i\in\Bbb Z}\g\langle i\rangle  \ ,
\]
where $\g\langle i\rangle =\{\,x\in\g\mid [h,x]=ix\,\}$.
Set $\g\langle{\ge}j\rangle=\oplus_{i\ge j}\g\langle i\rangle $.
We also write $\el$ for $\g\langle 0\rangle$ and $\p$ for
$\g\langle{\ge} 0\rangle$.
Since all characteristics of $e$ are $Z_G(e)$-conjugate,
the properties of this $\Bbb Z$-grading
do not depend on a particular choice of $h$.
\\[.5ex]
The orbit $G{\cdot}h$ contains a unique element $h_+$ such that
$h_+\in\frak t$ and $\alpha(h_+)\ge0$
for all $\alpha\in\Pi$.
The Dynkin diagram of $\g$ equipped with the numerical
labels $\alpha_i(h_+)$, $\ap_i\in\Pi$, at the corresponding nodes
is called the
{\it weighted Dynkin diagram} of $e$. After Dynkin, it is known
(see  \cite[8.1,\,8.3]{dynkin}) that \par
(a) $\ap_i(h_+)\in\{0,1,2\}$; \par
(b) $\tri$-triples $\{e,h,f\}$ and $\{e',h',f'\}$
are $G$-conjugate {\sl if and only if\/}
$h$ and $h'$ are $G$-conjugate {\sl if and only if\/}
their weighted Dynkin diagrams coincide.
\\[.5ex]
The following facts on the structure of $Z_G(e)\subset G$ and
$\z_\g(e)\subset\g$ are standard, see \cite[ch.\,III]{ss} or \cite{CoMc}.
\begin{s}{Proposition}         \label{stab}
Let $L$ (resp. $P$) be the connected subgroup of $G$ with Lie algebra
$\el$ (resp. $\p$), and put $K = Z_{G}(e)\cap L$. Then
\begin{itemize}
\item[{\sf (i)}]  $K=Z_G(e)\cap Z_G(f)=Z_G(f)\cap L$, and it is a maximal
reductive subgroup in both $Z_G(e)$ and $Z_G(f)$; \
$Z_G(e)\subset P$;
\item[{\sf (ii)}] the Lie algebra $\z_\g(e)$ (resp.\ $\z_\g(f)$) is
positively (resp. negatively) graded; $\z_\g(e)=
\bigoplus_{i\ge 0} \z_\g(e)\langle i\rangle $,
where $\z_\g(e)\langle i\rangle =\z_\g(e)
\cap\g\langle i\rangle $, and likewise for $\z_\g(f)$;
\item[{\sf (iii)}] for any $i$, there are $K$-stable
decompositions:
\[ \g\langle i\rangle =\z_\g(e)\langle i\rangle
   \oplus[f,\g\langle i+2\rangle\,]\qquad
   \g\langle i\rangle =\z_\g(f)\langle i\rangle  \oplus[e,\g\langle i-2\rangle
   \,].
\]
In particular, $\ad e:\g\langle i-2\rangle\to\g \langle i\rangle $ is injective
for $i\le 1$ and surjective for $i\ge 1$;
\item[{\sf (iv)}] $(\ad e)^i : \g\langle -i\rangle \rar\g\langle i\rangle $ is one-to-one;
\item[{\sf (v)}] $\dim\zge=\dim \g\langle 0\rangle+\dim \g\langle 1\rangle $.
\end{itemize}
\end{s}%
It follows that the reductive group $K$ is the centraliser of the
$\tri$-triple $\{e,h,f\}$.
Notice that $K$ can be disconnected and
$\zge\langle 0\rangle=\z_\g(f)\langle 0\rangle=\ka=\Lie K$.
By \cite[1.2]{p99},\\[.6ex]
\centerline{$\g\langle i\rangle $ is an orthogonal (resp. symplectic)
$K$-module if $i$ is even (resp. odd).}
In particular, $\dim\g\langle i\rangle $ is even
whenever $i$ is odd.
Let $S$ be a generic stabiliser for the orthogonal
$K$-module $\g\langle 2\rangle $. By a result of Luna, it is again a
reductive group. Formulas for the complexity and rank of $G{\cdot}e$
(see \cite[2.3]{p94} or \cite[4.2]{p99}) exploit
the above $\Bbb Z$-grading and the groups $K$ and $S$. We recall
it in Section~\ref{h3}.
\\[.5ex]
The integer $\max\{i\in {\Bbb N}\mid \g\langle i\rangle \ne 0\}$ is called the
{\it height\/} of $e$ or of the orbit $G{\cdot}e$ and is denoted by $\hot(e)$.
It is clear that $\hot(e)\ge 2$ for any $e\ne 0$, and
$\hot(e)=m$ if and only if $(\ad e)^m\ne 0$ and
$(\ad e)^{m+1}= 0$. See
\cite[Sect.\,2]{p99} for some results concerning the height.
\\[.6ex]
Whenever we discuss in the sequel
a $\Bbb Z$-grading associated with some $e\in\N$, this means that $e$ is
regarded as member of an $\tri$-triple and the grading is determined by the
corresponding semisimple element.

\begin{s}{Proposition} \label{odd-open}
Suppose $\hot(e)$ is odd, say $2d+1$. Then $K$ has an open orbit
in $\g\langle 2d{+}1\rangle $.
\end{s}\begin{proof} By a result of Vinberg \cite[2.6]{vi76},
$L$ has finitely many orbits in each $\g\langle i\rangle $,  $i\ne 0$.
Consider the periodic grading of $\g$ 
that is obtained by assembling together the spaces $\g\langle i\rangle $ modulo
$2d+3$. It is formally defined
by the inner automorphism $\vartheta=\Ad (\exp(\frac{2\pi\sqrt{-1}}{2d+3}h))$.
Let $\zeta$ be a primitive root of unity, of degree $2d+3$.
Then $\g=\displaystyle\bigoplus_{i=0}^{2d+2}\g_i$, where $\g_i$ is the
eigenspace of $\vartheta$ corresponding to the eigenvalue $\zeta^i$.
For a suitable choice of $\zeta$, we have $\g_0=\g\langle 0\rangle$, $\g_1=\g\langle 1\rangle $, and
$\g_2=\g\langle 2\rangle \oplus\g\langle -2d{-}1\rangle$. Therefore $L=G_0$ has an open orbit in $\g_1$.
It follows from Vinberg's argument (see \cite[3.1]{vi76}) that
the same holds for any $\g_i$ such that
$\gcd(i,2d+3)=1$. In particular,
$L$ has an open orbit in $\g_2$. By the definition of $K$, it is
equivalent to the fact that $K$ has an open orbit in $\g\langle -2d{-}1\rangle$.
Since $K$ is reductive, the same holds for $\g\langle 2d{+}1\rangle $.
\end{proof}%
The symbol $\odot$ is used below for almost direct products of
algebraic groups. This means that one has a direct sum for the corresponding
Lie algebras.

\begin{s}{Proposition}   \label{2d+1}
Suppose $\hot(e)=2d+1$. Then
\begin{itemize}
\item[{\sf 1}.] $\g\langle 2d{+}1\rangle =\oplus_{i=1}^s \VV_i$ and
$K=\prod_{i=1}^s Sp(\VV_i) \odot K_1$, where the reductive group $K_1$
acts trivially on $\g\langle 2d{+}1\rangle $.
\item[{\sf 2}.]
$L=GL(\g\langle 2d{+}1\rangle )\odot L_1$, where the reductive group $L_1$
acts trivially on $\g\langle 2d{+}1\rangle $. In particular,
$\g\langle 2d{+}1\rangle$ is a two-orbit
$L$-module.
\end{itemize}
\end{s}\begin{proof*} 1. Since $\g\langle 2d{+}1\rangle$ is a symplectic $K$-module,
this is a straightforward consequence of  Theorem~\ref{main1} and
Proposition~\ref{odd-open}.  \\
2. Let $\co$ be the open $K$-orbit in $\g\langle 2d{+}1\rangle $. It follows from part 1 that
its complement has exactly $s$ irreducible components:
\begin{equation} \label{complement}
 \g\langle 2d{+}1\rangle \setminus\co=\bigcup_{j=1}^s (\bigoplus_{i\ne j}\VV_i)  .
\end{equation}
Let $\tilde\co\supset\co$ be the open $L$-orbit in $\g\langle 2d{+}1\rangle $.
By Eq.~\re{complement}, its complement is contained in a union of
{\it proper\/} subspaces. Because $\g$ is a simple Lie algebra,
$\g\langle 2d{+}1\rangle $ is a simple $L$-module. Therefore
this complement must be the origin (for, any non-trivial irreducible
component would generate a proper $L$-stable subspace).
Thus, $\g\langle 2d{+}1\rangle $ is a two-orbit $L$-module.

It is clear that the $\Bbb Z$-graded Lie subalgebra
$\g\langle -2d{-}1\rangle \oplus\g\langle 0\rangle\oplus\g\langle 2d{+}1\rangle $ is reductive. Set $\q=
[\,\g\langle 2d{+}1\rangle , \g\langle -2d{-}1\rangle \,]\subset \g\langle 0\rangle$. It is easily seen that $\h:=
\g\langle -2d{-}1\rangle \oplus\q\oplus\g\langle 2d{+}1\rangle $ is a (reductive) Lie algebra. Furthermore,
since the representation $\varrho: \g\langle 0\rangle \to {\frak gl}(\g\langle 2d{+}1\rangle )$ is
irreducible, $\h$ is simple and $\g\langle 0\rangle=\q\oplus \Ker(\varrho)$.
Let $Q\subset L$ be the connected group with Lie algebra $\q$.
Our goal is to prove that $Q\simeq GL(\g\langle 2d{+}1\rangle )$. Since $Q$ acts
transitively on $\g\langle 2d{+}1\rangle \setminus\{0\}$, this quickly follows
from the classification of 3-term gradings of simple Lie algebras.
An {\sf a priori} proof is given in the next Lemma.
\end{proof*}%
\begin{s}{Lemma}  \label{glv}
Let $\h=\h\langle -1\rangle\oplus\h\langle 0\rangle\oplus\h\langle 1\rangle $ be a simple $\Bbb Z$-graded Lie algebra.
Suppose the group $H\langle 0\rangle$ acts transitively on $\h\langle 1\rangle \setminus\{0\}$.
Then $H\langle 0\rangle\simeq GL(\h\langle 1\rangle )$ and $H\simeq GL(V')$, where $\dim V'=\dim\h\langle 1\rangle +1$.
\end{s}\begin{proof}
Set $V=\h\langle 1\rangle $. Then $\h\langle -1\rangle\simeq V^*$.
Let $\beta$ (resp. $\theta$) be the lowest (resp. highest) weight of
the $H\langle 0\rangle$-module $V$. (We assume that $\te_H\subset \h\langle 0\rangle$ and
$\n_+\supset \h\langle 1\rangle =V$.) It is clear that $\theta$ is the highest root and
$\beta$ is a simple root of $\h$. Let $\Pi'$ be the set of all other simple
roots. Write $\theta=n_\beta\beta+\sum_{\gamma\in\Pi'}n_\gamma\gamma$.
It is well known that $\beta$ is long and $n_\beta=1$.
\\
By assumption, $[\h,v]=V$ for any $v\in V\setminus\{0\}$. It follows that
$[\xi,v]\ne 0$ for any $v\in V\setminus\{0\}$ and $\xi\in V^*\setminus\{0\}$.
In particular, $[\xi_{-\beta}, v_\theta]\ne 0$. Hence $\theta-\beta$ is
a (positive) root of $\h\langle 0\rangle$ and $(\theta,\beta^\vee)=(\theta^\vee,\beta)
=1$. We have
\[
  2=(\theta,\theta^\vee)=1+ (\sum_{\gamma\in\Pi'}n_\gamma\gamma,
  \theta^\vee) \ .
\]
It follows that there exists a unique $\beta'\in\Pi'$ such that
$(\beta',\theta^\vee)=1$ and $n_{\beta'}=1$. This means that $\beta'$
is long as well.
It is well known that the support of any root on the
Dynkin diagram is connected. Therefore $\beta$ and $\beta'$ are
extreme roots on the Dynkin diagram of $\h$. Hence
$\theta-\beta$ is the highest root in the {\sl irreducible\/} root system with
basis $\Pi'$. Let $\beta_2$ be the unique
simple root adjacent to $\beta$ on the Dynkin diagram. If $\beta_2=\beta'$,
then $\h\simeq {\frak sl}_3$, and we are done. Otherwise, we obtain
$((\theta-\beta)^\vee,\beta_2)=((\theta-\beta)^\vee,\beta')=1$, so that we may
argue by induction on the length of the unique chain connecting $\beta$ and
$\beta'$ on the Dynkin diagram. Finally, we obtain that
$\theta=\sum_{1\le i\le n}\beta_i$, where $\beta=\beta_1$,
$\beta'=\beta_n$, and $\Pi'=\{\beta_2,\dots,\beta_n\}$.
From this, all assertions follows easily.
\end{proof}%
{\bf Remark.} The integer $s$ appearing in Proposition~\ref{2d+1} is the
number of irreducible constituents of the $K$-module $\g\langle 2d{+}1\rangle $.
Explicit description of the nilpotent elements with odd height shows
that one always has $s=1$, i.e., $\g\langle 2d{+}1\rangle $ is a simple $K$-module.
In case $d=1$, this can be proved
{\sf a priori}, see Section~\ref{h3}.

\section{On nilpotent orbits of height 3}
\label{h3}
\setcounter{equation}{0}
\noindent
Recall that, for any irreducible $G$-variety $X$, one defines the complexity
of $X$, denoted $c_G(X)$, and the rank of $X$, denoted $r_G(X)$.
Here $c_G(X)$ is the minimal
codimension of orbits in $X$ of a Borel subgroup of $G$. If $X$ is quasiaffine
(which we only need), then $r_G(X)$ is the dimension of the $\bbk$-vector
space in $\te^*_{{\Bbb Q}}$ generated
by the highest weights of all simple $G$-modules occurring in $\bbk[X]$, see
\cite[ch.\,1]{diss} for all this.
\\[.5ex]
Keep the notation introduced in Section~\ref{z-graded}.
In the context of nilpotent orbits, the following was proved in
\cite[Theorem~2.3]{p94}:
\begin{equation} \label{c}
 c_G(G{\cdot}e)=c_L(L/K) + c_S(\g\langle{\ge}3\rangle  ) \ ,
\end{equation}
\begin{equation} \label{r}
 r_G(G{\cdot}e)=r_L(L/K) + r_S(\g\langle{\ge}3\rangle  ) \ .
\end{equation}
Let $G{\cdot}e$ be an orbit of height 3. Here $\ka$ is a symmetric
subalgebra of $\el$ \cite[3.3]{p94};
hence $L/K$ is an $L$-spherical homogeneous space. In this particular case,
Eq.~\re{c} says that
\begin{itemize}
\item[]  The orbit $G{\cdot}e$ is $G$-spherical if and only if
$\g\langle 3\rangle$ is a spherical $S$-module, i.e., a Borel subgroup of $S$ has an
open orbit in $\g\langle 3\rangle$.
\end{itemize}

\begin{s}{Theorem} \label{height3}
Suppose $\hot(e)=3$. Then $G{\cdot}e$ is spherical.
\end{s}\begin{proof}
By Proposition~\ref{2d+1}, we have $\g\langle 3\rangle=\oplus_{i=1}^s \VV_i$,
$K=\prod_{i=1}^s Sp(\VV_i) \odot K_1$,
and $L=GL(\g\langle 3\rangle) \odot L_1$. Here $K_1\subset L_1$
and the composition
\vskip.6ex
\centerline{
$\prod_{i=1}^s Sp(\VV_i)\hookrightarrow L\to GL(\oplus_{i=1}^s \VV_i)$}
\vskip.7ex\noindent
yields the standard embedding.
Because $L/K$ is a spherical homogeneous space, we see that
$GL(\oplus_{i=1}^s \VV_i)/\prod_{i=1}^s Sp(\VV_i)$ is also a spherical
homogeneous space. Since the dimension of the space of
$\prod_{i=1}^s Sp(\VV_i)$-fixed vectors in the $GL(\g\langle 3\rangle)$-module
$\wedge^2\g\langle 3\rangle$ equals $s$, one must have
$s=1$, i.e.,
$K=Sp(\g\langle 3\rangle)\odot K_1$. Moreover,
since $\ka$ is symmetric in $\el$, the group $Sp(\g\langle 3\rangle)$ cannot be
embedded diagonally in $GL(\g\langle 3\rangle)\odot L_1$. Thus,
$Sp(\g\langle 3\rangle)\subset GL(\g\langle 3\rangle)$ and $K_1\subset L_1$.
In particular, $\ka_1$ is symmetric in $\el_1$.
Then the $K$-module $\g\langle 2\rangle $ is isomorphic
to $\Lie L/\Lie K=\g\langle 0\rangle/\ka\simeq
\wedge^2(\g\langle 3\rangle)\oplus (\el_1/\ka_1)$, where
$Sp(\g\langle 3\rangle)$ acts only on $\wedge^2(\g\langle 3\rangle)$ and
$K_1$ act only on $\el_1/\ka_1$.
Recall that a generic stabiliser of the $Sp(\VV)$-module
$\wedge^2\VV$ is isomorphic to $(SL_2)^{\dim\VV /2}$.
Therefore $S\simeq (SL_2)^d \times S_1$, where $d=\dim\g\langle 3\rangle/2$ and
$S_1$ is a generic stabiliser for the $K_1$-module
$\el_1/\ka_1$.
Since $S_1$ acts trivially on $\g\langle 3\rangle$,
the explicit structure of this group is unimportant.
Thus, $\g\langle 3\rangle$ is isomorphic as $S$-module to the
sum of simplest representations of different $SL_2$.
Now it is clear that $\g\langle 3\rangle$ is a spherical $S$-module, and we are done.
\end{proof}%
Using the notation of the preceding proof, we can give a formula for the rank
of $G{\cdot}e$.
\begin{s}{Proposition}  \label{rank}
If\/ $\hot(e)\le 3$, then $r_G(G{\cdot}e)=r_{L_1}(L_1/K_1)+\dim\g\langle 3\rangle$.
\end{s}\begin{proof} We apply Eq.~\re{r}.
Since $\g\langle 3\rangle$ as $S$-module is a sum of simplest representations of different
$SL_2$, we have $r_S(\g\langle 3\rangle)=(\dim\g\langle 3\rangle)/2$.
Since $L/K\simeq GL(\g\langle 3\rangle)/Sp(\g\langle 3\rangle) \times L_1/K_1$, we obtain
$r_L(L/K)=(\dim\g\langle 3\rangle)/2 + r_{L_1}(L_1/K_1)$.
\end{proof}%
{\it Remark.} If $\hot(e)=2$, then $\g\langle 3\rangle=0$, $L_1=L$, $K_1=K$, and
the formula applies as well. Since $L_1/K_1$ is a symmetric space, its rank
is a familiar object. E.g., the rank is immediately seen from the Satake
diagram.
\\[.6ex]
Let $Q$ be a connected algebraic group with Lie algebra $\q$.
Consider the coadjoint representation of $Q$ in $\q^*$.
The transcendence degree of the field $\bbk(\q^*)^Q$ is called
the {\it index\/} of $\q$. Recall the following two basic results
on the index (see \cite{p02} for more details).

$\bullet$ -- \un{Ra\"\i s' formula for the index of semi-direct products}:
Let $V$ be a $Q$-module.
The vector space sum  $\q\oplus V$ has a natural
structure of a Lie algebra such that $V$ is an Abelian ideal and
$\q$ is a subalgebra; the bracket of $\q$ and $V$ is determined by the
representation of $\q$ on $V$. The resulting semi-direct product
is denoted by $\q\, \semidir\, V$.
Then
\[
\ind(\q\,\semidir\, V)=\ind\q_\nu + \trdeg \bbk(V^*)^Q \,
\]
where $\nu\in V^*$ is a generic point.

$\bullet$ -- \un{Vinberg's inequality}: if $\xi\in \q^*$ is arbitrary, then
$\ind\q_\xi\ge \ind\q$.

\begin{s}{Theorem} \label{index3}
Suppose $\hot(e)=3$. Then $\ind\zge=\rk\g$.
\end{s}\begin{proof}
Our plan is as follows. Since $\g\simeq\g^*$, Vinberg's inequality says
that $\ind\zge\ge\rk\g$. We find a special point $\xi\in\zge^*$ such that
$\zge_\xi$ has a semi-direct product structure. Using Ra\"\i s' formula
and the above structure results for nilpotent orbits of height 3,
we compute that $\ind\zge_\xi=\rk\g$. Then the second application of
Vinberg's inequality (with $\q=\zge$) shows that $\ind\zge\le \rk\g$.
\\
The centraliser of $e$ has the following graded structure:
\[
\zge=\ka\oplus\zge\langle 1\rangle \oplus\g\langle 2\rangle \oplus\g\langle 3\rangle \ .
\]
Using $\Phi$, consider the opposite element of the $\tri$-triple, $f$, as a
linear form on $\zge$. It is an easy exercise based on Proposition~\ref{stab}
that $\zge_f=\ka\oplus \g\langle 2\rangle \oplus\g\langle 3\rangle$. In other words,
$\zge_f\simeq \ka\,\semidir\,(\g\langle 2\rangle \oplus\g\langle 3\rangle)$.
Write $V$ for $\g\langle 3\rangle$. We know that
$\g\langle 0\rangle\simeq \glv\oplus\el_1$ and
$\ka\simeq \spv\oplus \ka_1$, where $\el_1$ acts trivially on $V$ and
$\ka_1\subset\el_1$.
Moreover, $\ka_1$ is a symmetric subalgebra of $\el_1$ and
$\g\langle 2\rangle \simeq \wedge^2 V \oplus (\el_1/\ka_1)$ as $\ka$-module.
\\[.6ex]
The above description shows that $\zge_f$ is a direct sum of two
semi-direct products:
\[
  \zge_f= (\spv\,\semidir\,(\wedge^2 V\oplus V))\dotplus
  (\ka_1\,\semidir\,(\el_1/\ka_1)) \ ,
\]
where `$\dotplus$' stands for the direct sum of Lie algebras.
By Ra\"\i s' formula, the index of the first summand is $\dim V$
(a generic stabiliser of $\spv$ in $\wedge^2 V^*\oplus V^*$
is a {\sl commutative\/} subalgebra of dimension $(\dim V)/2$; hence the
transcendence degree of the field $\bbk(\wedge^2 V^*\oplus V^*)^{Sp(V)}$ is
also $(\dim V)/2$).
The index for the second summand is computed in the proof of Theorem 3.5 in
\cite{p02}. The only important thing here is that
$\ka_1\to {\frak so}(\el_1/\ka_1)$ is the isotropy representation of the
symmetric subalgebra $\ka_1\subset\el_1$.
This implies that the second index is equal to $\rk\el_1$.
Finally, the above formula for $\g\langle 0\rangle$ shows that $\dim V+\rk\el_1=\rk\g\langle 0\rangle=
\rk\g$.
\end{proof}%
Together with \cite[3.5]{p02}, this Theorem shows that $\ind\zge=\rk\g$
whenever $G{\cdot}e$ is spherical.

\section{On a specific nilpotent orbit of height 3}
\label{special}
\setcounter{equation}{0}
\noindent
It was noticed in \cite{p94} that orbits of height 3 exist in all simple Lie
algebras except $\slv$ and $\spv$. Furthermore, the height of any nilpotent
orbit in $\slv$ or $\spv$ is even (see \cite[2.3(1)]{p99}).
Here we give a partial explanation of this phenomenon by giving an intrinsic
construction of a
nilpotent orbit of height 3 in each simple Lie algebra whose
highest root is fundamental, i.e., for
$\g\ne\spv$ or $\slv$. This nilpotent orbit, which will be denoted by
$\Bbb O$, satisfies a number of
interesting relations, see below.
\\[.6ex]
For $\gamma\in \Delta$, let $e_\gamma\in\g$ be a nonzero root vector.
We assume that $\Phi(e_\gamma,e_{-\gamma})=1$. If $\gamma$ is long, then
the orbit $G{\cdot}e_\gamma$ is the minimal nonzero nilpotent orbit,
but the word `nonzero' is usually omitted in this case. Consider the
$\tri$-triple $\{e_\theta, h_\theta=[e_\theta, e_{-\theta}], e_{-\theta}\}$
and the corresponding $\Bbb Z$-grading
$\displaystyle\bigoplus_{-2\le i\le 2}\g\langle i\rangle _{(\theta)}$. As before, we write
$L$ for the connected group with Lie algebra $\g\langle 0\rangle_{(\theta)}$, etc.
This grading was considered by many authors and for various purposes.
Our point of view is invariant theoretic, and
we want to deduce our results in a classification-free way.
Here $\g\langle 2\rangle _{(\theta)}=\bbk e_\theta$ and therefore $\ka$ is a subalgebra of
codimension 1 in $\el$. More precisely, $\el=\ka\oplus \bbk h_\theta$.
From now on, we assume that $\theta$ is a fundamental weight.
It is worth mentioning that, assuming this constraint and also that
$\rk\g\ge 4$,
G.\,R\"ohrle \cite{gerhard} describes natural bijections between the set of
$L$-orbits in $\gthe$ and certain double cosets of $G$ and its Weyl group.

\begin{s}{Proposition}   \label{for_fundamental}
Suppose $\theta$ is a fundamental weight. Then
\begin{itemize}
\item[{\sf (i)}] $\g\langle 1\rangle _{(\theta)}$ is a simple symplectic $K$-module;
\item[{\sf (ii)}] the quotient space $\g\langle 1\rangle _{(\theta)}\md K$ is 1-dimensional;
\item[{\sf (iii)}] the algebra $\bbk[\g\langle 1\rangle _{(\theta)}]^K$ is generated
by a single polynomial of degree 4, say $F$. This polynomial is
explicitly given by
$F(x)=\Phi((\ad x)^4 e_{-\theta}, e_{-\theta})$, $x\in\g\langle 1\rangle _{(\theta)}$.
\end{itemize}
\end{s}\begin{proof}
(i) We only need to prove that $\g\langle 1\rangle _{(\theta)}$ is simple.
For any root vector $e_\nu\in\g$, we have $[h_\theta,e_\nu]=
(\nu,\theta^\vee)e_\nu$.
Let $\beta$ be the {\sl unique\/} simple root such that $(\theta,\beta)> 0$.
Then $(\theta^\vee,\beta)=1$, i.e., $e_\beta\in \g\langle 1\rangle _{(\theta)}$.
Since $[\ka,\ka]=[\el,\el]$, it is clear that $e_\beta$ is the unique lowest
weight vector in $\g\langle 1\rangle _{(\theta)}$ with respect to $\ka\cap \ut_+$.
\\[.6ex]
(ii) As $L$ has finitely many orbits in $\g\langle 1\rangle _{(\theta)}$, we see that
$\dim (\g\langle 1\rangle _{(\theta)}\md K)\le 1$. Clearly, $\theta-\beta$ is the highest weight of the $L$-module
$\g\langle 1\rangle _{(\theta)}$. For $\nu\in\te^*$, let $\bar\nu$ denote the restriction of
$\nu$ to $\te\cap\ka$. Then $\ov{\theta-\beta}$ (resp. $\bar\beta$)
is the highest (resp. lowest) weight of the $\ka$-module $\g\langle 1\rangle _{(\theta)}$.
(Notice that $\ov{\theta}=0$.)
Because $\theta$ is the fundamental weight
corresponding to $\beta$, we have $(\beta^\vee,\theta)=1$.
Hence $(\beta^\vee,\theta)=(\theta^\vee,\beta)$ and therefore
$\beta$ is a long root.
Hence $\theta-2\beta\not\in\Delta$, i.e. $\ov{\theta-\beta}-\ov{\beta}$
is not a root of $\ka$. It then follows from Corollary~\ref{w} that
$\bbk[\g\langle 1\rangle _{(\theta)}]^K\ne \bbk$.
\\[.6ex]
(iii) It follows from (ii) that  $\bbk[\g\langle 1\rangle _{(\theta)}]^K$ is a polynomial
algebra. Since $\g\langle 1\rangle _{(\theta)}$ is simple and symplectic, there are no
$K$-invariants of degree 1 and 2. On the other hand,
Proposition~\ref{degree-4} says that there is an invariant of degree 4.
Therefore, this invariant must generate $\bbk[\g\langle 1\rangle _{(\theta)}]^K$.
A formula for $F$ in terms of $\Phi$ and a quadratic mapping $\psi$ is given
in the proof of \re{degree-4}.
In our case,  $\psi: \g\langle 1\rangle _{(\theta)} \to\ka\subset \g\langle 0\rangle_{(\theta)}$
is given by $\psi(x)=(\ad x)^2 e_{-\theta}$. Then the required
formula for $F$ follows by the $\g$-invariance of $\Phi$.
\end{proof}%
{\bf Remarks.} 1. If $\theta$ is fundamental,
then $K=(L,L)$, and it is connected and semisimple.
\par
2. For $\g=\spn$, we have $\theta=2\vp_1$. Here $\bbk[\g\langle 1\rangle _{(\theta)}]^K=\bbk$.
\par
3. For $\g=\sln$ ($n\ge 2$), we have $\theta=\vp_1+\vp_n$.
Here the $K$-module $\g\langle 1\rangle _{(\theta)}$ is not simple, and
$\bbk[\g\langle 1\rangle _{(\theta)}]^K$ is generated by a polynomial of degree 2.

\begin{s}{Theorem}  \label{postroenie}
Suppose $\theta$ is fundamental. Then
\begin{itemize}
\item[{\sf (i)}] The open $L$-orbit in $\g\langle 1\rangle _{(\theta)}$, say $\tilde\co$,
is affine. We have $\tilde\co=\{x\in\gthe\mid \hot(x)=4\}$.
If $\tilde e$ is any element of $\tilde\co$, then $2h_\theta$ is a
characteristic of $\tilde e$,
$\dim(\Ima(\ad \tilde e)^4)=1$, and $\dim(\Ima(\ad \tilde e)^3)=2$.
\item[{\sf (ii)}] The complement of $\tilde\co$ in $\gthe$ is an irreducible
variety. If $\co$ is the dense $L$-orbit in
$\g\langle 1\rangle _{(\theta)}\setminus\tilde\co$, then $\hot(e)=3$ and
$\dim(\Ima(\ad e)^3)=2$ for any $e\in\co$.
\end{itemize}
\end{s}\begin{proof}
(i) Since $K=(L,L)$ has invariants on $\gthe$, we can apply
\cite[Prop.\,3.3(1)]{kac}. That Proposition describes how $L$-orbits break
into $K$-orbits.
Namely, the open $L$-orbit is affine and it consists of a
1-parameter family of {\it closed\/} $K$-orbits. All other $K$-orbits are
unstable and they coincide with $L$-orbits. Thus,
the complement of the open $L$-orbit in $\g\langle 1\rangle _{(\theta)}$ is given by
the equation $F=0$. From Proposition~\ref{for_fundamental}(iii), it follows
that $F(x)\ne 0$ if and only if $(\ad x)^4e_{-\theta}\ne 0$, i.e.,
$\hot(x)=4$. For $\tilde e\in\tilde\co$, we have $\z_\el(\tilde e)=
\z_\ka(\tilde e)$. It is a reductive Lie algebra, since $K{\cdot}\tilde e$
is closed. Therefore $\z_\el(\tilde e)$
is orthogonal to $h_\theta$ with respect to $\Phi$. Then Proposition~1.2 in
\cite{kac} says that the line $\bbk h_\theta$ is in the image of $\ad\tilde e$.
Hence $2h_\theta$ is a characteristic of $\tilde e$. The last two equalities
easily follow from this and the fact that $\dim\g\langle 2\rangle _{(\theta)}=1$.
\par
(ii) Since $K$ is semisimple, $F$ is an irreducible polynomial. Hence
$\gthe\setminus\tilde\co$ is irreducible.
If $F(x)=0$, then $\hot(x)\le 3$. However it is not yet clear that
there does exist a point with height 3. To see this,
consider the differential of $F$. It is easy to derive a formula for $dF_x$.
If $y\in\gthe$, then
\begin{equation}  \label{differential}
dF_x(y)= 2\,\Phi\Bigl((\ad y)(\ad x)^3e_{-\theta},e_{-\theta}\Bigr)+
2\,\Phi\Bigl((\ad x)(\ad y)(\ad x)^2 e_{-\theta},e_{-\theta}\Bigr) \ .
\end{equation}
If $\hot(x)=2$, then the first summand is zero. Some standard transformations
show that the same is true for the second summand. Indeed,
$(\ad x)(\ad y)(\ad x)^2 e_{-\theta}=(\ad y)(\ad x)^3e_{-\theta}+
(\ad [x,y])(\ad x)^2 e_{-\theta}$. Notice that $[x,y]=c e_\theta$. If
$c=0$, then we are done. If not, then up to a scalar factor,
the second summand in Eq.~\re{differential}
is  equal to
\[
\Phi((\ad e_\theta)(\ad x)^2e_{-\theta},e_{-\theta})=
-\Phi((\ad x)^2e_{-\theta}, h_\theta)=-\Phi((\ad x)e_{-\theta}, x)=0 .
\]
Thus, $\hot(e)=3$ whenever  $dF_e\ne 0$. The latter holds for $e$ in the dense
$L$-orbit in the hypersurface $\{F=0\}$.
Since $G{\cdot}e$ lies in the closure of $G{\cdot}\tilde e$, we have
$0< \dim(\Ima(\ad e)^3)\le \dim(\Ima(\ad\tilde e)^3)=2$.
Consider a $\Bbb Z$-grading associated with $e$: \ $\g=\displaystyle
\bigoplus_{-3\le i\le 3}\g\langle i\rangle $.
Since $\dim\g\langle 3\rangle$ is even (see Sect.~\ref{z-graded}) and
$\g\langle 3\rangle=\Ima(\ad e)^3$, the dimension in question is 2.
\end{proof}%
Set $\tilde{\Bbb O}=G{\cdot}\tilde\co$
and ${\Bbb O}=G{\cdot}\co$.
We have characterised $\tilde{\Bbb O}$ as the nilpotent orbit whose
weighted Dynkin diagram is twice that
of the minimal nilpotent orbit. We also described
${\Bbb O}$ as the unique nilpotent orbit that is dense in
$\gthe\setminus (\tilde{\Bbb O}\cap\gthe)$.
Below, we give a more direct characterisation of ${\Bbb O}$.

\begin{s}{Proposition}  \label{unique}
\begin{itemize}
\item[1.] \ $G{\cdot}\gthe$ is closed in $\g$;
\item[2.] \ For any $x\in \gthe\setminus\{0\}$, we have $\dim G{\cdot}x=
2\dim L{\cdot}x+2$;
\item[3.] \ ${\Bbb O}$ is the unique $G$-orbit of codimension 2 in the closure
of ${\tilde{\Bbb O}}$.
\end{itemize}
\end{s}\begin{proof}
1. Since $\g\langle{\ge}1\rangle_{(\theta)}$ is the nilpotent radical of a
parabolic, $G{\cdot(\g\langle{\ge}1\rangle_{(\theta)})}$ is closed.
Recall that $\g\langle{\ge}1\rangle_{(\theta)}$ is a Heisenberg Lie algebra (see \cite{jo})
and $\g\langle 2\rangle _{(\theta)}=\bbk e_\theta$.
Therefore if $x=x_1+x_2\in\g\langle{\ge}1\rangle$ and $x_1\ne 0$, then $G{\cdot}x$ has a
representative in $\gthe$. Finally, $G{\cdot}e_\theta$ has a representative in
$\gthe$, since $\gthe$ contains long root spaces.
\par
2. Let us look at the graded structure of $\z_\g(x)$.
Using the Kostant-Kirillov form associated with $x$, one sees that
$\dim\g\langle i\rangle_{(\theta)}-\dim\z_\g(x)\langle i\rangle =
\dim\g\langle -i{-}1\rangle_{(\theta)}-\dim\z_\g(x)\langle -i{-}1\rangle$
for all $i$.
Clearly, $\z_\g(x)\langle 2\rangle=\g\langle 2\rangle _{(\theta)}$, $\dim\z_\g(x)\langle 1\rangle =\dim\gthe-1$,
and $\dim L{\cdot}x=\dim\g\langle 0\rangle_{(\theta)}-\dim\z_\g(x)\langle 0\rangle$. Now, the conclusion
is obtained by a simple counting.
\par
3. By part (1), $\ov{\tilde{\Bbb O}}=G{\cdot}\gthe$, and we have proved
in Theorem~\ref{postroenie}(ii)  that
$\co$ is the unique $L$-orbit of codimension 1 in $\gthe$.
\end{proof}%
{\it Remark.}
Since $\Bbb O$ is $G$-spherical, Proposition~\ref{unique}(3) implies that
$c_G(\tilde{\Bbb O})\le 2$. It can conceptually be proved that the
complexity in question equals 2. Moreover, $\tilde{\Bbb O}$
is {\it the\/} unique minimal non-spherical orbit, if $\theta$ is
fundamental and $\g\ne\GR{F}{4}$, see \cite[4.4]{p99}.  \\
Now, we are able to describe the weighted Dynkin diagram
and other properties of ${\Bbb O}$.

\begin{s}{Theorem}   \label{char_3}
Let $e\in {\Bbb O}$ and let $\displaystyle\bigoplus_{i=-3}^3 \g\langle i\rangle $ be the
$\Bbb Z$-grading associated with an $\tri$-triple $\{e,h,f\}$. Then
\begin{itemize}
\item[{\sf (i)}] $\dim\g\langle 3\rangle=2$ and $\dim\g\langle 1\rangle =2\dim\g\langle 2\rangle $;
\item[{\sf (ii)}]
Let $\beta\in\Pi$ be the unique root such that $(\theta,\beta)>0$, and
let $\{\ap_i\}$ ($i\in I$) be all simple roots adjacent to $\beta$ on the
Dynkin diagram. Then weighted Dynkin diagram of $\Bbb O$ is obtained by
attaching `1'  to the all $\ap_i$'s $(i\in I)$ and `0' to all other simple
roots.
\end{itemize}
\end{s}\begin{proof*}
(i) The first equality is already proved in Theorem~\ref{postroenie}(ii).
To prove the second, we construct a bi-grading of $\g$.
Take a nonzero element $x\in\g\langle 3\rangle$. By a modification of the Morozov-Jacobson
theorem, one can find an $\tri$-triple $\{x,\tilde h,y\}$ such that
$\tilde h\in \g\langle 0\rangle$ and $y\in\g\langle -3\rangle$.
Recall that $\g\langle 3\rangle\setminus \{0\}$ is a single $L$-orbit, see
Proposition~\ref{2d+1}(2). Therefore $x$ lies in the minimal nilpotent
orbit and $\hot(x)=2$. The pair of commuting elements
$h,\tilde h$ determines a refinement of the original $\Bbb Z$-grading:
\[
        \g=\bigoplus_{i,j}\g\langle i,j\rangle  \ ,
\]
where $\g\langle i,j\rangle =\{z\in\g\mid [h,z]=iz,\ [\tilde h,z]=jz \}$,
$-3\le i\le 3$, and $-2\le j\le 2$. We write $\g\langle i,\ast\rangle $ for
$i$-th eigenspace of $h$ and $\g\langle \ast,j\rangle $ for $j$-eigenspace
of $\tilde h$. Notice that $x\in\g\langle 3,2\rangle $ and therefore $\ad x$
takes $\g\langle i,j\rangle $ to $\g\langle i+3,j+2\rangle$.
We are going to show that many subspaces of this bi-grading
are equal to zero.

By Proposition~\ref{stab}(iv), the linear mappings
$(\ad x)^2: \g\langle \ast,-2\rangle \to \g\langle \ast,2\rangle $ and $\ad x: \g\langle \ast,-1\rangle \to \g\langle \ast,1\rangle $
are one-to-one; therefore $\g\langle i,-2\rangle =0$ for $i\ge -2$,
$\g\langle i,-1\rangle =0$ for $i\ge 1$, and $\dim\g\langle i,-1\rangle
=\dim\g\langle i+3,1\rangle $ for $i\le 0$.
We also have $[\g\langle 0,\ast\rangle ,x]=\g\langle 3,\ast\rangle $, hence
$\g\langle 3,0\rangle=\{0\}$.

Using the central symmetry for dimensions and the fact that
$\dim\g\langle 3,\ast\rangle =2$, we obtain
\[
  \dim\g\langle 3,2\rangle =\dim\g\langle 3,1\rangle =\dim\g\langle 0,-1\rangle =\dim\g\langle 0,1\rangle =1
\]
and
\[
  \dim\g\langle 2,1\rangle =\dim\g\langle -1,-1\rangle =\dim\g\langle 1,1\rangle =\dim\g\langle -2,-1\rangle =:a
\]
Thus, we obtain the following ``matrix" of numbers
$\dim\g\langle i,j\rangle $, where only possible non-zero entries are indicated.
\begin{center}
\begin{tabular}{r|rrrrrrr}
2     &    &     &     &     &     &     & 1  \\
1     &    &     &     & 1   & $a$ & $a$ & 1  \\
0     &    & $c$ & $b$ & $d$ & $b$ & $c$ &    \\
-1    & 1  & $a$ & $a$ & 1   &     &     &    \\
-2    & 1  &     &     &     &     &     &    \\ \hline
$j/i$ & \phantom{\hspace{3ex}}-3   & \phantom{\hspace{3ex}}-2
      & \phantom{\hspace{3ex}}-1   & \phantom{\hspace{3ex}}0
      & \phantom{\hspace{3ex}}1    &  \phantom{\hspace{3ex}}2
      & \phantom{\hspace{3ex}}3  \\
\end{tabular}
\end{center}
Since $e$ lies in $\g\langle 2,0\rangle \oplus\g\langle 2,1\rangle $ and
$\ad e:\g\langle -1,\ast\rangle \to\g\langle 1,\ast\rangle $ is one-to-one,
we see that $\ad e$ maps $\g\langle -1,-1\rangle $ injectively to
$\g\langle 1,0\rangle $. Hence $a\le b$. Now we compute $\dim\zge$ in two
different ways. Using $i$-grading and Proposition~\ref{stab}(v), we obtain
$\dim\zge=\dim\g\langle 0,\ast\rangle +\dim\g\langle 1,\ast\rangle =2+a+b+d$.
On the other hand, the  $j$-grading arises
as grading connected with the minimal nilpotent orbit. Therefore
$2\tilde h$ is a characteristic
for $\tilde{\Bbb O}$ (Theorem~\ref{postroenie}(i))
and hence $\dim\z_\g(\tilde e)=\dim\g\langle \ast,0\rangle =
2b+2c+d$ for any $\tilde e\in \tilde{\Bbb O}$. Because
$\dim{\Bbb O}=\dim\tilde{\Bbb O}-2$ (Proposition~\ref{unique}(3)),
we obtain $2b+2c+d=a+b+d$. Hence $c=0$ and $a=b$.
\par
(ii)
Let $y$ be a nonzero element of the 1-dimensional space
$\g\langle 0,1\rangle $. Obviously, the ``highest"
subspace in each column lies in $\z_\g(y)$; next,
a subspace of codimension 1 in $\g\langle 0,0\rangle $ also lies there. Hence
$\dim\z_\g(y)\ge 4a+d+2$. On the other hand, we deduce from considering
the $j$-grading that $\dim\z_\g(x)=4a+d+2$. Therefore $y$ lies in the minimal
nilpotent orbit as well, and one has equality for $\dim\z_\g(y)$.
This implies that $\ad y:\g\langle 1,0\rangle \to\g\langle 1,1\rangle $ is one-to-one.

Without loss of generality, we may assume that the triangular decomposition
of $\g$ is chosen so that $h,\tilde h\in \te\subset \g\langle 0,0\rangle $ and
$\ut_+\supset \g\langle i,j\rangle $ for all $(i,j)$ such that $i+j>0$.
Then $x$ becomes a highest weight
vector with respect to this choice of $\Pi$, i.e.,
$x=e_\theta$.
Since $y$ is a lowest weight vector in the {\it simple\/}
$\g\langle \ast,0\rangle $-module $\g\langle \ast,1\rangle $, neither of the root
spaces in $\g\langle i,1\rangle $, $i\ge 1$, corresponds to simple roots.
Thus, the simple root spaces
lie only in $\g\langle 0,0\rangle ,\,\g\langle 0,1\rangle $ and $\g\langle 1,0\rangle $. In particular, $y=e_\beta$
for some $\beta\in\Pi$. Let $\{\ap_i\}$ be all simple roots such that
the corresponding root spaces belong to $\g\langle 1,0\rangle $.
The argument in the previous paragraph and bi-graded structure of $\g$
show that $\beta$
is the only simple root such that $(\beta,\theta)\ne 0$ and the
$\ap_i$'s are exactly the simple roots that are not orthogonal to
$\beta$. Thus, $\ap(h)\le 1$ for all $\ap\in\Pi$, and $\ap(h)=1$
if and only if $(\ap,\beta)\ne 0$.
\end{proof*}%
\begin{rem}{Remarks}  \label{apparent} 1.
Once the equalities $c=0$, $b=a$ are proved, one may observe that the
above ``dimension matrix" gains the apparent $\GR{G}{2}$-symmetry. For,
after making an affine transformation, the matrix can be depicted as
follows:
\begin{center}
\begin{tabular}{r|rcccccc}
$j$   &    &     &     & \\ \hline
2     &    &     &     &  1  &     &     &    \\
1     &  \phantom{\hspace{3ex}}1 &     & $a$ &     & $a$ &     & 1  \\
0     &    & $a$ &     & $d$ &     & $a$ &    \\
-1    &  \phantom{\hspace{3ex}}1 &     & $a$ &     & $a$ &     & 1  \\
-2    &    &     &     &  1  &     &     &    \\
\end{tabular}
\end{center}
\vskip1ex\noindent
Each one-dimensional space represents the highest root space with respect to
a suitable choice of $\Delta^+$. It is not hard to prove that 6 one-dimensional
spaces generate a Lie subalgebra of type $\GR{A}{2}$. Next, our nilpotent
element $e$ of height 3 lies in one of $a$-dimensional spaces. It can be shown
that 6 one-dimensional spaces together with $e$ generate a Lie algebra of type
$\GR{G}{2}$. Therefore the ${\Bbb Z}\times {\Bbb Z}$-grading constructed
in the previous proof yields also an instance of a Lie algebra
graded by root system $\GR{G}{2}$. (See \cite{josiane} and references
therein for the general notion of a Lie algebra graded by a root system.)
In this way, one obtains an intrinsic construction
of a $\GR{G}{2}$-grading for any
simple Lie algebra whose highest root is fundamental.
\end{rem}%
2. Explicit calculations show that one always has
$r_G({\Bbb O})=r_G(\tilde{\Bbb O})=\min\{ \rk\g , 4\}$.

\section{Regular functions on a nilpotent orbit}
\setcounter{equation}{0}
\noindent
In this section, we propose an approach to describing the algebra of
covariants on nilpotent orbits.
\\[.7ex]
Maintain notation of Section~\ref{z-graded}.
Let $\{e,h,f\}$ be an arbitrary $\tri$-triple. Without loss, we may assume that
$\ap(h)\ge 0$ for all $\ap\in\Pi$ (i.e., $h=h_+$). Then
$\ut_+=\g\langle 0\rangle_+\oplus \g\langle {\ge} 1\rangle$ and
$\ut_-=\g\langle 0\rangle_-\oplus \g\langle {\le}{-}1\rangle$,
where $\g\langle 0\rangle_\pm =\ut_\pm \cap \g\langle 0\rangle$.
Let $\Delta\langle i\rangle$ denote the subset of $\Delta$ corresponding to
$\g\langle i\rangle$.
Set
\[
U_-=\exp(\ut_-),\ U=\exp(\ut_+),\ U(L)_-=\exp(\g\langle 0\rangle_-),\
U(L)=\exp(\g\langle 0\rangle_+) \ .
\]
These are maximal unipotent subgroups of $G$ and $L$, respectively.
Then $U(L)_-=U_-\cap P=U_-\cap L$ and
$U_-\simeq N\times U(L)_-$, where $N=\exp(\g\langle {\le}{-}1\rangle)$.
Since $\gedva$ is a $P$-module, one can
form the homogeneous vector bundle $G\ast_P\gedva$ over $G/P$. The following is
fairly well known:
\begin{equation} \label{funkcii}
     \bbk[ G\ast_P\gedva]=\bbk[G{\cdot}e]=\bbk[(\ov{G{\cdot}e})_n] ,
\end{equation}
where $(\ )_n$ denotes the normalisation, and
\begin{equation} \label{collaps}
   \textrm{the collapsing }\quad G\ast_P\gedva \to \ov{G{\cdot}e} \quad
\textrm{ is a projective birational morphism}.
\end{equation}
It follows from Eq.~\re{funkcii} and \re{collaps} that
$\gedva$ can be regarded as closed subvariety of $(\ov{G{\cdot}e})_n$.
Therefore there is the restriction homomorphism
$\bbk[G{\cdot}e] \to \bbk[\gedva]$, which is onto. Clearly, it takes
$U_-$-invariant functions on $G{\cdot}e$ to $U(L)_-$-invariant functions
on $\gedva$. Thus, we obtain the restriction homomorphism
\[
   \tau^o : \bbk[G{\cdot}e]^{U_-} \to \bbk[\gedva]^{U(L)_-} \ .
\]
Notice that $T=\exp(\te)$ is a maximal torus in both $L$ and $G$, and that
both subalgebras under consideration are $T$-stable.

\begin{s}{Theorem} \label{restrict}
$\tau^o$ is $T$-equivariant, injective and birational.
\end{s}\begin{proof}
The first claim is obvious.
Since $U_-{\cdot}P$ is dense in $G$ and hence $U_-{\cdot}\gedva$ contains
a dense open subset of $G{\cdot}e$, we see that $\tau^o$ is injective.
To prove birationality,
we dualize the picture and consider the morphism of the corresponding spectra:
\[
   \tau :  \gedva\md U(L)_- \to G{\cdot}e\md U_- \ .
\]
The birationality of $\tau$ is equivalent to the fact that there exists
a dense open subset $\mathcal D$ of $\gedva$ such that \\
\hbox to \textwidth{\quad $(\ast)$\hfil
$U_-{\cdot}v\cap \gedva=U(L)_-{\cdot}v$ for all $v\in {\mathcal D}$.
\hfil}
We take ${\mathcal D}=P{\cdot}e$. To prove $(\ast)$,
we first notice that \
$(\ast\ast)$\quad$G{\cdot}e\cap\gedva=P{\cdot}e$. Indeed,
$P{\cdot}e$ is dense in $\gedva$ in view of Proposition~\ref{stab}(iii),
and if $v\in\gedva\setminus P{\cdot}e$, then $\dim Z_P(v) >\dim Z_P(e)=
\dim Z_G(e)$.
Second, we notice that if $v\in P{\cdot}e$, $g\in N$, and $g{\cdot}v\in
\gedva$, then $g=1$. This follows from $(\ast\ast)$ and the containment
$Z_G(e)\subset P$. Now, $(\ast)$ easily follows.
\end{proof}%
The following observation will allow us to switch between
$\bbk[G{\cdot}e]^U$ and $\bbk[G{\cdot}e]^{U_-}$.

\begin{s}{Lemma}   \label{dual}
Given $e\in \N$, let $\Gamma=\Gamma(G{\cdot}e)$ be the monoid of all
highest weights occurring in $\bbk[G{\cdot}e]$, i.e. the monoid of $T$-weights
in $\bbk[G{\cdot}e]^{U}$.
Then the monoid of all lowest weights in $\bbk[G{\cdot}e]$ is $-\Gamma$.
\end{s}\begin{proof}
By Frobenius reciprocity, this is equivalent to the fact that if $V$ is a
simple $G$-module and $V^{Z_G(e)}\ne 0$, then
$(V^*)^{Z_G(e)}\ne 0$. The last assertion follows from the following
two properties of a Weyl involution $\vartheta\in\Aut G$:

(a) the $G$-module $V$ equipped with the twisted action  $(g,v)\mapsto
\vartheta(g){\cdot}v$, $g\in G$, $v\in V$, is isomorphic to $V^*$;

(b) $\vartheta$ acts trivially on the set of nilpotent orbits \cite[2.10]{ls}.
In particular, $\vartheta(Z_G(e))$ is $G$-conjugate to $Z_G(e)$.
\end{proof}%
Given a (locally-finite) $T$-module $\ca$, let $\ca_\mu$ denote the
$\mu$-weight space of $\ca$. We have $\displaystyle \bbk[G{\cdot}e]^{U_-}=
\bigoplus_{\mu\in -\Gamma}\bbk[G{\cdot}e]^{U_-}_\mu$ and
$\tau^o: \bbk[G{\cdot}e]^{U_-}_\mu\hookrightarrow \bbk[\gedva]^{U(L)_-}_\mu$.
\\
Let ${\mathcal X}_+(G)$ (resp. ${\mathcal X}_-(G)$) denote the cone of
dominant (resp. antidominant) weights of $G$.
The mapping $\tau^o$ is not always onto.
There are two related difficulties:

(1) since $L$ has fewer simple roots than $G$, the cone of (anti)dominant
weights for $L$ is wider than that for $G$. So that $L$-lowest weight vectors
corresponding to ``extra-weights" (if any) cannot lie in the image of $\tau^o$;

(2) even for $\mu\in {\mathcal X}_-(G)$,
the dimensions of $\bbk[G{\cdot}e]^{U_-}_\mu$ and
$\bbk[\gedva]^{U(L)_-}_\mu$ can be different.
\\
To overcome the first difficulty, one may consider the subalgebra
of $\bbk[\gedva]^{U(L)_-}$ corresponding to ${\mathcal X}_-(G)$.
Still, this does not yield an isomorphism in general
(e.g., in case of the principal nilpotent
orbit in ${\frak sl}_3$), because of the presence of the second problem.
If $G{\cdot}e$ is spherical, then $\gedva$ is a spherical $L$-module
\cite[4.2(2)]{p99}. Hence all {\it nonzero\/} weight spaces in
$\bbk[G{\cdot}e]^{U_-}$ and $\bbk[\gedva]^{U(L)_-}$ are 1-dimensional.
This partly resolves the second problem and gives some hope
that a more precise statement holds in the spherical case.
To begin with, consider the orbits of height 2.
Then $\ov{G{\cdot}e}$ is normal by a result of W.~Hesselink \cite{he79}.
Therefore there is no difference
between $\bbk[G{\cdot}e]$ and $\bbk[\ov{G{\cdot}e}]$. Hesselink used Kempf's
theorem on the collapsing of homogeneous bundles with completely reducible
fibres. Since this result applies not only to nilpotent orbits in
$\g$, we will work for a while in a more general setting.
\\[.7ex]
Let $P$ be a parabolic subgroup of $G$ and let $U_-$ be a maximal unipotent
subgroup of $G$ such that $U_-{\cdot}P$ is dense in $G$.
Then $P\cap U_-=:U(L)_-$ is a maximal unipotent subgroup for some Levi subgroup
$L\subset P$.
\begin{s}{Lemma}   \label{kempf}
Let $V$ be a $G$-module and $M\subset V$ a $P$-submodule. Suppose
$M$ is completely reducible and the collapsing
$\pi: G\ast_PM \to G{\cdot}M\subset V$ is birational. Then \\[.6ex]
\centerline{
$\bbk[G{\cdot}M]^{U_-}\simeq \bbk[G\ast_PM]^{U_-}\simeq
\bbk[M]^{U(L)_-}$ \ .}
\vskip-.8ex
\end{s}\begin{proof}
Since $\pi$ is proper and birational, we have
$\bbk[G{\cdot}M]\simeq \bbk[G\ast_PM]$,
and hence $\bbk[G{\cdot}M]^{U_-}\simeq \bbk[G\ast_PM]^{U_-}$.
\\[.6ex]
The second map is given by the restriction to the fibre over $\{P\}$
of the projection $G\ast_PM \to G/P$. Obviously, this restriction is
injective (even if $M$ is not completely reducible). Ontoness follows by
an application of the Borel-Weil-Bott theorem that is due to G.\,Kempf
(see Hesselink's version in \cite[sect.\,3]{he79}). Namely,
\[
\bbk[G\ast_PM]\simeq \oplus_{n\ge 0}H^0(G/P, {\mathcal L}(S^nM^*)) \ .
\]
As $S^nM^*$ is a completely reducible $P$-module and $S^nV^* \to
S^nM^*$ is onto, each simple $P$-module (= $L$-module) $E$ in
$S^nM^*$ gives rise  to a nontrivial simple $G$-module
$H^0(G/P, {\mathcal L}(E))$ in $\bbk[G\ast_PM]$.
\end{proof}%
Now, we come back to nilpotent orbits.

\begin{s}{Theorem}  \label{height2}
Suppose $\hot(e)=2$ and hence $\gedva=\g\langle 2\rangle$. Then
\begin{itemize}
  \item[{\sf (i)}]\quad $\tau^o: \bbk[\ov{G{\cdot}e}]^{U_-}\isom
\bbk[\g\langle 2\rangle]^{U(L)_-}$, and it is a polynomial algebra;
  \item[{\sf (ii)}]\quad $\bbk[\g\langle 2\rangle]$ is a free
$\bbk[\g\langle 2\rangle]^{U(L)}$-module;
  \item[{\sf (iii)}]\quad $\bbk[\ov{G{\cdot}e}]$ is a free
$\bbk[\ov{G{\cdot}e}]^{U}$-module.
\end{itemize}
\end{s}\begin{proof*}
(i) Since $G{\cdot}e$ is spherical, $\g\langle 2\rangle$ is a spherical
$L$-module.
That $\bbk[\g\langle 2\rangle]^{U(L)_-}$ is a polynomial algebra is a
standard property of spherical representations. Since $\g\langle 2\rangle$
is a completely reducible $P$-module, ontoness of $\tau^o$
follows from Lemma~\ref{kempf}.
\par
(ii) This is essentially proved in \cite[4.6]{aura}. To get in that
situation, one has to consider the ``even'' reductive subalgebra
$\g^{ev}=\g\langle -2\rangle\oplus \g\langle 0\rangle\oplus
\g\langle 2\rangle\subset\g$. Then $\g\langle 2\rangle$ is the Abelian
nilpotent radical of $\p^{ev}=\g\langle 0\rangle\oplus \g\langle 2\rangle$.
\par
(iii) In \cite[5.5]{aura}, we proved a general sufficient condition for
$\bbk[X]$ to be a flat $\bbk[X]^U$-module, where $G$ is reductive and
$X$ is an affine $G$-variety. This reads as follows: \\[.6ex]
\hspace*{2pt}$(\spadesuit)$\quad \parbox{435pt}{%
Suppose $\bbk[X]^U$ is polynomial and let $\lb_1,\dots,\lb_r$ be the
$T$-weights of free generators of $\bbk[X]^U$. Suppose
the $\lb_i$'s are linearly independent and any fundamental
weight of $G$ occurs in at most one $\lb_i$ (i.e., the $\lb_i$'s depend on
pairwise disjoint subsets of fundamental weights). Then $\bbk[X]$ is a flat
$\bbk[X]^U$-module. (If $G$ is semisimple, then the linear independence of
the $\lb_i$'s follows from the second assumption.)
}
\\[.6ex]
In the graded situation, ``free'' is equivalent to ``flat''.
Hence our aim is to verify this property for $X=\ov{G{\cdot}e}$. Using the
isomorphism in part (i) and some results of \cite{aura},
one can give the explicit description of the $\lb_i$'s.
Let $\gamma_1,\dots,\gamma_r$ be the {\it upper canonical string\/}
(=\,u.c.s.) of roots in $\Delta\langle 2\rangle$.
By definition, this means that $\gamma_1=\theta$; $\gamma_2$ is {\it the\/}
maximal root in $\Delta^{(1)}\langle 2\rangle=\{\ap\in
\Delta\langle 2\rangle\mid (\ap,\gamma_1)=0\}$, and so on... The procedure
terminates when $\Delta^{(i)}\langle 2\rangle$ becomes empty. See
\cite[Sect.\,1]{aura} for more details.
(Our $\Delta\langle 2\rangle$ here is $\Delta(1)$ in \cite{aura}. In course of
constructing the u.c.s., one may ignore the subset
$\Delta\langle 1\rangle$.)
This string consists of pairwise orthogonal long positive roots.
\\
It follows from \cite[Corollary~4.2]{aura} that the weights of free generators
of $\bbk[\g\langle 2\rangle]^{U(L)_-}$ are $-\gamma_1,-\gamma_1-\gamma_2,
\dots,-\gamma_1-\cdots-\gamma_r$. More precisely, that Corollary describes
the weights for $\bbk[\g\langle 2\rangle]^{U(L)}$ in terms of the
{\it lower canonical string\/} (=\,l.c.s.). Using the longest element in the
Weyl group of $L$, it is then easy to realise that
the description of $\bbk[\g\langle 2\rangle]^{U(L)_-}$ can be given in terms
of the u.c.s., as above.
Then using (i) and Lemma~\ref{dual},
we conclude that
the weights of free generators of $\bbk[\ov{G{\cdot}e}]^U$,
i.e., the generators of $\Gamma$, are
\begin{equation}  \label{generators}
\lb_1=\gamma_1,\lb_2=\gamma_1+\gamma_2,\dots,
\lb_r=\gamma_1+\cdots+\gamma_r \ .
\end{equation}
Thus, we are to
verify that, for any $\ap\in\Pi$, there is at most one $i\in\{1,\dots,r\}$
such that $(\ap,\lb_i)> 0$. The proof of this
will be given in the next Lemma.
\end{proof*}%
\begin{s}{Lemma}   \label{ucs}
Suppose $\hot(e)=2$ and
let $\gamma_1,\dots,\gamma_r$ be the u.c.s.
in $\Delta\langle 2\rangle$. Then for any $\ap\in\Pi$,
there is at most one $i\in\{1,\dots,r\}$
such that $(\ap,\gamma_1+\dots+\gamma_i)> 0$.
\end{s}\begin{proof*}
Let $\Pi\langle i\rangle$ denote the subset of $\Pi$ which lies in
$\Delta\langle i\rangle$. Then $\Pi\langle 0\rangle$ is the set of simple
roots for $L$. For $\ap\in\Pi\langle 0\rangle$, the desired property
has been proved in \cite[4.6]{aura}. And this is the key property that
yields the freeness in \cite{aura} and hence in Theorem~\ref{height2}(ii).
Dealing with the other simple roots, one has
to distinguish two possibilities.

{\bf (a)} \ $\g\langle 1\rangle=0$, i.e. $e$ is an even nilpotent element.
\\
Then there is a unique root $\beta\in\Pi\langle 2\rangle$.
This $\beta$ is the lowest weight of the $L$-module $\g\langle 2\rangle$
with respect to $U(L)$.
Since the $\Bbb Z$-grading of $\g$ is
associated with of a nilpotent orbit, one can show that the u.c.s. and l.c.s.
in $\Delta\langle 2\rangle$ coincide as sets (cf. also Prop.\,1.5
in \cite{aura}).
In particular, this means that $\beta=\gamma_r$. Hence $(\beta,\lb_i)=0$ if
$i<r$ and $(\beta,\lb_r)=(\beta,\beta)>0$.

{\bf (b)} \ $\g\langle 1\rangle\ne 0$.
\\
Then there are one or two roots in $\Pi\langle 1\rangle$.
Let $\ap\in \Pi\langle 1\rangle$.
It is clear that $(\ap,\gamma_i)\ge 0$ for all $i$, since $\Delta\langle
j\rangle=\varnothing$ for $j\ge 3$.
\par
$(b_1)$ There is at most one index $i$ such that $(\ap,\gamma_i)>0$.
\\
Indeed, if $(\ap,\gamma_i)>0$ and $(\ap,\gamma_j)>0$, then $\gamma_i-\ap
\in\Delta\langle 1\rangle$ and $(\gamma_i-\ap,\gamma_j)<0$. Whence
$\gamma_i+\gamma_j-\ap\in \Delta\langle 3\rangle$. A contradiction!
\par
$(b_2)$ Assume that $(\ap,\gamma_i)>0$ for $i<r$. Then we also have
$\gamma_{i+1}$. Since $\gamma_{i+1}\ne \theta$, there is $\nu\in\Pi
\langle 0\rangle$ such that $\gamma_{i+1}+\nu$ is a root (in $\Delta
\langle 2\rangle$). Notice that $\gamma_{i+1}+\nu\ne \gamma_i$, since
$\gamma_i$ and $\gamma_{i+1}$ are orthogonal long roots. By the definition of
the u.c.s., we have $(\gamma_i,\gamma_{i+1}+\nu)\ne 0$ and hence it is
positive. Therefore $(\gamma_i,\nu)>0$ and \\
\hbox to \textwidth{ \
$(\Diamond)$ \hfil $\gamma_i-\gamma_{i+1}-\nu\in\Delta\langle 0\rangle^+$.
\hfil}
\\[.7ex]
(This root is positive, since $\hot(\gamma_i)> \hot(\gamma_{i+1})$.)
We have $\gamma_i-\ap\in\Delta\langle 1\rangle$ and
$(\gamma_i-\ap,\gamma_{i+1}+\nu)=(\gamma_i,\nu)-(\ap,\nu)>0$. Hence \\
\hbox to \textwidth{ \
$(\Diamond\Diamond)$ \hfil
$\gamma_{i+1}+\nu-\gamma_i+\ap\in\Delta\langle 1\rangle$.
\hfil}
\\[.7ex]
Taking the sum of $(\Diamond)$ and $(\Diamond\Diamond)$, we conclude that
$\ap$ is not simple. A contradiction!
Thus, $(\ap,\gamma_i)=0$ if $i<r$ and $(\ap,\gamma_r)\ge 0$,
which implies the assertion.
\end{proof*}%
\begin{s}{Complement {\ququ (to Theorem~\ref{height2}(i))}} \label{dop}
Let $f_i\in \bbk[\ov{G{\cdot}e}]^{U}$ be the generator with weight $\lb_i$.
Then $\deg f_i=i$.
\end{s}\begin{proof}
The isomorphism $\tau^o$ preserve the degree of functions, and the degrees
of generators of $\bbk[\g\langle 2\rangle]^{U(L)}$ were described in
\cite[4.1]{aura}.
\end{proof}%
{\bf Remark.}
For $G=SL_n$, the fact that $\bbk[G{\cdot}e]^U$ is polynomial was previously
proved in \cite[Sect.\,8]{mitya}. Shmelkin also explicitly describes the
generators as certain minors of $n$ by $n$ matrices.
\\[.7ex]
Making use of Theorem~\ref{height2}, Eq.~\re{generators} and
Complement~\ref{dop}, it is quite
easy to determine the structure of $\bbk[\ov{G{\cdot}e}]$ as graded $G$-module,
if $\hot(e)=2$. The relevant information will be presented in Table~1 below.
\par
Set
$\displaystyle\bbk[\gedva]^{U(L)_-}_{(-)}=\bigoplus_{\mu\in {\mathcal X}_-(G)}
\bbk[\gedva]^{U(L)_-}_\mu$. Clearly, $\tau^o
(\bbk[G{\cdot}e]^{U_-})\subset \displaystyle\bbk[\gedva]^{U(L)_-}_{(-)}$.

\begin{s}{Conjecture}   \label{conj1}
Suppose $G{\cdot}e$ is spherical. Then
$\tau^o:
\bbk[G{\cdot}e]^{U_-}\isom \displaystyle\bbk[\gedva]^{U(L)_-}_{(-)}$.
\end{s}%
The preceding exposition shows that Conjecture~\ref{conj1} is true
if $\hot(e)=2$, and one even has $\bbk[\g\langle 2\rangle]^{U(L)_-}_{(-)}=
\bbk[\g\langle 2\rangle]^{U(L)_-}$.
But for non-spherical orbits the above map can fail to be an isomorphism.
Thus, the only interesting open case is that of orbits of height~3.
It is more convenient to work with $U$-invariants, because then the weights
of generators are dominant. To this end, Conjecture~\ref{conj1} can be restated
in the ``dominant'' form:
\begin{s}{Conjecture$'$}   \label{conj1prim}
Suppose $G{\cdot}e$ is spherical. Then
$\hat\tau^o:
\bbk[G{\cdot}e]^{U}\isom \bbk[\g\langle {\le}{-}2\rangle]^{U(L)}_{(+)}
\simeq {\cs}(\gedva)^{U(L)}_{(+)}$.
\end{s}%
Here $\cs(\cdot)$ stands for the symmetric algebra, and the subscript $(+)$
means that we take only weight spaces from ${\mathcal X}_+(G)$.
\\[.6ex]
We say that a monoid $\Gamma\subset {\mathcal X}_+(G)$ is
{\it saturated\/}, if ${\Bbb Z}\Gamma\cap {\mathcal X}_+(G)=\Gamma$.
For instance, it follows from Eq.~\re{generators} and
Lemma~\ref{ucs} that $\Gamma(G{\cdot}e)$ is saturated if $\hot(e)=2$.

\begin{s}{Conjecture}   \label{conj2}
If $e\in\N$ is arbitrary, then $\Gamma(G{\cdot}e)$ is saturated.
\end{s}%
{\it Remark.} If $G{\cdot}e$ is the regular nilpotent orbit, then
$\bbk[G{\cdot}e]=\bbk[\N]$. Here $\Gamma(\N)$ is intersection
of the root lattice with ${\mathcal X}_+(G)$; hence it is
saturated. More generally, for $e$ even, it follows from \cite{mac} that
$\Gamma(G{\cdot}e)=\Gamma(G/L)$. This means it is interesting to
study saturatedness property for $G/L$ with an arbitrary Levi subgroup $L$.
It is likely that Littelmann's path method can be helpful for the last
problem.
\\[.7ex]
Using methods of \cite{diss},
it can be shown that, for the spherical orbits, Conjecture~\ref{conj2}
implies Conjecture~\ref{conj1}. This would provide an easy way for
determining the structure of the algebra of covariants and, in particular,
the monoid $\Gamma(G{\cdot}e)$ for {\sl all\/} spherical orbits, see
Section~\ref{tables}.

\section{Tables and examples}  \label{tables}
\setcounter{equation}{0}
\noindent
In Tables 1 and 2, we collect some information on spherical nilpotent
orbits.
Numbering of the simple roots and fundamental weights of $\g$
is the same as in \cite[Tables]{vion}.
In the classical case, we follow the standard notation for roots
via $\esi_i$'s. Roots of exceptional Lie algebras
are presented via the Dynkin diagrams. For instance,
the array
$\left(\!\!\begin{array}{l} n_1\ n_2\ n_3\ n_4\ n_5\ n_6\ n_7 \\
              \phantom{n_1\ n_2\ n_3\ n_4\ }n_8  \end{array}
      \!\!\right)$
represents the root $\sum n_i\ap_i$ for $\GR{E}{8}$. This picture also
demonstrates our numbering of fundamental weights for $\GR{E}{8}$.
In both Tables, nilpotent orbits are described by partitions (the classical
case) and weighted Dynkin diagrams (the exceptional case).
\\[.7ex]\indent
In Table 1, we present the information on the monoid
$\Gamma(G{\cdot}e)$ for the nilpotent orbits of height 2.
The relations between the $\gamma_i$'s and the $\lb_i$'s are given
in Eq.~\re{generators}.
As is well-known, a very even partition in the orthogonal case
gives rise to two nilpotent orbits. This happens if $l=0$ for the second
item in the ${\frak so}_n$-case. However, we include in the Table only
one possibility. The second possibility is (only for $l=0$):
$\gamma_r=\esi_{2r-1}-\esi_{2r}$ and $\lb_r=2\vp_{2r-1}$.
We also did not include in Table~1 the minimal
nilpotent orbits for the exceptional Lie algebras, because we only have
$\gamma_1=\lb_1=\theta$ in these cases.

\begin{table}[htb]
\begin{center}
\caption{ 
Nilpotent orbits of height 2}
\begin{tabular}{|c|lll|}  \hline
$\g$ & nilpotent orbit  & \quad u.c.s. & Free generators of $\Gamma$ \\
\hline \hline
${\frak sl}_{n}$ & $(2^r,1^{n-2r})$
      & $\gamma_i=\esi_i-\esi_{n-i+1}$   & $\lb_i=\vp_i+\vp_{n-i}$       \\
  &  $r\le n/2$    & \ \ $i=1,\dots,r$ &    \\ \hline
${\frak sp}_{2n}$ & $(2^r,1^{2n-2r})$
      & $\gamma_i=2\esi_i$     & $\lb_i=2\vp_i$       \\
  &  $r\le {n}$    & \ \ $i=1,\dots,r$ &  \\  \hline
${\frak so}_{n}$ & $(3,1^{l})$
      & $\gamma_1=\esi_{1}+\esi_{2}$   & $\lb_1=\vp_{2}$    \\
  &  $l\ge 4$, $n=3{+}l$
      & $\gamma_2=\esi_1-\esi_{2}$     & $\lb_2=2\vp_1$   \\ \cline{2-4}
  &   $(2^{2r},1^{l})$
      & $\gamma_i=\esi_{2i-1}+\esi_{2i}$   & $\lb_i=\vp_{2i}$ $(i<r)$  \\
  & \raisebox{12pt}{$n=4r{+}l$}
      & \raisebox{12pt}{\ \ $i=1,\dots,r$}     &
   $\lb_r{=}\left\{\!\begin{array}{ll} 2\vp_{2r}, & l\le 1 \\
                                 \vp_{2r}{+}\vp_{2r+1},   &  l=2   \\
                                   \vp_{2r},  &  l\ge 3
                 \end{array}\right.$
    \\ \hline
$\GR{E}{6}$ & {\small\begin{tabular}{@{}c@{}}
1--0--\lower3.5ex\vbox{\hbox{0\rule{0ex}{2.5ex}}
\hbox{\hspace{0.4ex}\rule{.1ex}{1ex}\rule{0ex}{1.4ex}}\hbox{0\strut}}--0--1
\end{tabular}}  &
     $\gamma_1=\genfrac{(}{)}{0pt}{}{1\,2\,3\,2\,1}{2}$
                  &  $\lb_1=\vp_6$    \\
     &  &
     $\gamma_2=\genfrac{(}{)}{0pt}{}{1\,1\,1\,1\,1}{0}$
                  &  $\lb_2=\vp_1+\vp_5$   \\[1ex] \hline
$\GR{E}{7}$ & {\small\begin{tabular}{@{}c@{}}
0--1--0--\lower3.5ex\vbox{\hbox{0\rule{0ex}{2.5ex}}
\hbox{\hspace{0.4ex}\rule{.1ex}{1ex}\rule{0ex}{1.4ex}}\hbox{0\strut}}--0--0
\end{tabular}}  &
      $\gamma_1=\genfrac{(}{)}{0pt}{}{1\,2\,3\,4\,3\,2}{\phantom{w}2}$
                &  $\lb_1=\vp_6$    \\
     &  &
      $\gamma_2=\genfrac{(}{)}{0pt}{}{1\,2\,2\,2\,1\,0}{\phantom{w}1}$
                &  $\lb_2=\vp_2$    \\[1ex] \cline{2-4}
      & {\small\begin{tabular}{@{}c@{}}
2--0--0--\lower3.5ex\vbox{\hbox{0\rule{0ex}{2.5ex}}
\hbox{\hspace{0.4ex}\rule{.1ex}{1ex}\rule{0ex}{1.4ex}}\hbox{0\strut}}--0--0
\end{tabular}}  &
    $\gamma_1=\genfrac{(}{)}{0pt}{}{1\,2\,3\,4\,3\,2}{\phantom{w}2}$
        &  $\lb_1=\vp_6$    \\
     &  &
       $\gamma_2=\genfrac{(}{)}{0pt}{}{1\,2\,2\,2\,1\,0}{\phantom{w}1}$
                   &  $\lb_2=\vp_2$    \\[1ex]
     &  &
       $\gamma_3=\genfrac{(}{)}{0pt}{}{1\,0\,0\,0\,0\,0}{\phantom{w}0}$
                  &  $\lb_3=2\vp_1$    \\[1ex] \hline
$\GR{E}{8}$  & {\small\begin{tabular}{@{}c@{}}
0--0--0--0--\lower3.5ex\vbox{\hbox{0\rule{0ex}{2.5ex}}
\hbox{\hspace{0.4ex}\rule{.1ex}{1ex}\rule{0ex}{1.4ex}}\hbox{0\strut}}--0--1
\end{tabular}}  &
  $\gamma_1=
   \genfrac{(}{)}{0pt}{}{2\,3\,4\,5\,6\,4\,2}{\phantom{\hspace*{2.1ex}}3}$
                      &  $\lb_1=\vp_1$    \\
     &  &
  $\gamma_2=
   \genfrac{(}{)}{0pt}{}{0\,1\,2\,3\,4\,3\,2}{\phantom{\hspace*{2.1ex}}2}$
                    &  $\lb_2=\vp_7$    \\[1ex]  \hline
$\GR{F}{4}$ & $1{-}0{\Leftarrow} 0{-}0$
      & $\gamma_1=(2\ 4\ 3\ 2)$   & $\lb_1=\vp_4$       \\
  &
      & $\gamma_2=(2\ 2\ 1\ 0)$   & $\lb_2=2\vp_1$   \\ \hline
\end{tabular}
\end{center}
\vskip-1ex
\end{table}
\vskip1ex 
In Table~2, we present our computations of the homogeneous generators of
$\cs(\g\langle 2\rangle\oplus\g\langle 3\rangle)^{U(L)}_{(+)}$
for all nilpotent orbits of height 3. By Conjecture~\ref{conj1prim},
these should also be the generators of the algebra $\bbk[G{\cdot}e]^U
=\bbk[(\ov{G{\cdot}e})_n]^U$.
Each generator is represented by its $T$-weight and degree, the latter
being the usual degree in the symmetric  algebra.
The corresponding monoid of dominant weights is called
$\tilde\Gamma$, and the third column gives the dimension of
the subspace in $\te^*_{\Bbb Q}$ generated by $\tilde\Gamma$ or,
equivalently, $r_G(G{\cdot}e)$ or the Krull dimension of $\bbk[G{\cdot}e]^U$.
Then using Theorem~\ref{char_3}(ii) or Remark~\ref{apparent}(2),
the interested reader may realise which orbit is $\Bbb O$ for each $\g$.
The rightmost column gives known information about
the normality of $\ov{G{\cdot}e}$. (The classical cases are due to
Kraft and Procesi; the case of $\GR{G}{2}$ is due to Levasseur and Smith;
the other exceptional cases are due to Broer, see \cite[p.\,959]{bram}.)
To fill in the fourth column, we first compute the weights and degrees of
the homogeneous generators of
$\cs(\g\langle 2\rangle\oplus\g\langle 3\rangle)^{U(L)}$. In doing so, one has
to take into account not only the fundamental weights of $L$, but also the
weights of the central torus in $L$. The main technical problem is then to
express the weights obtained in terms of fundamental weights of $G$. Once this
is done, it is rather easy to determine the monomials in  generators whose
weights are $G$-dominant.
\\[.6ex]
\begin{table}[htb]
\begin{center}
\caption{
The nilpotent orbits of height 3}
\begin{tabular}{|c|rcl|c|}  \hline
$\g$ & nilpotent orbit  & $\dim\tilde\Gamma$ & Weights and degrees of the
generators & Norm.\\
\hline \hline
$\begin{array}{c} {\frak so}_{4t+3} \\
                  (\GR{B}{2t+1})\end{array}$
            & $\begin{array}{r} (3,2^{2t}) \\
                 t\ge 1 \end{array}$ & $2t{+}1$ &
         $\begin{array}{ll}(\vp_{2i-1},i), & i=1,\dots,t \\
                  (\vp_{2i},i), & i=1,\dots,t \\
                  (2\vp_{2t+1}, t+1)
    \end{array}$ & $-$ \\ \hline
$\begin{array}{c} {\frak so}_{4t+4} \\
                  (\GR{D}{2t+2}) \end{array}$
            & $\begin{array}{r} (3,2^{2t},1)\\
                t\ge 1 \end{array}$ & $2t{+}2$ &
$\begin{array}{l}(\vp_{2i-1}+\vp_{2j-1},i+j),\ 1\le i\le j\le t \\
              (\vp_{2i},i),\ i=1,\dots,t \\
              (\vp_{2j-1}{+}\vp_{2t+1}{+}\vp_{2t+2},t{+}j{+}1),\ j=1,\dots,t \\
              (2\vp_{2t+1}, t+1), \ 
              (2\vp_{2t+2}, t+1)
    \end{array}$ & $+$ \\ \hline
${\frak so}_{4t{+}l{+}3}$ & $(3,2^{2t},1^l)$ & $2t{+}2$ &
                 \ $(\vp_{2i},i)$, $i=1,\dots,t$ & $+$\\
            & $t\ge 1,l\ge 2$ & &
            \ $(\vp_{2i-1}+\vp_{2j-1},i+j)$,  $1\le i\le j\le t{+}1$ & \\
       & & & $\left\{\!\begin{array}{rl}
                 (2\vp_{2t+2},t+1), & l=2 \\
                 (\vp_{2t+2}+\vp_{2t+3},t+1), & l=3 \\
                 (\vp_{2t+2},t+1), &  l\ge 4
         \end{array}\right.$ &   \\ \hline
$\GR{E}{6}$ &
{\small\begin{tabular}{@{}c@{}}
0--0--\lower3.5ex\vbox{\hbox{1\rule{0ex}{2.5ex}}
\hbox{\hspace{0.4ex}\rule{.1ex}{1ex}\rule{0ex}{1.4ex}}\hbox{0\strut}}--0--0
\end{tabular}}  & 4 & \ $(\vp_6,1),(\vp_1{+}\vp_5,2),(\vp_3,3),
               (\vp_2{+}\vp_4,4)$ & $+$ \\ \hline
$\GR{E}{7}$ & {\small\begin{tabular}{@{}c@{}}
0--0--0--\lower3.5ex\vbox{\hbox{0\rule{0ex}{2.5ex}}
\hbox{\hspace{0.4ex}\rule{.1ex}{1ex}\rule{0ex}{1.4ex}}\hbox{0\strut}}--1--0
\end{tabular}}  & 4 & \ $(\vp_6,1),(\vp_2,2),(\vp_5,3),
               (\vp_4,4)$ & + \\ \cline{2-5}
       & {\small\begin{tabular}{@{}c@{}}
1--0--0--\lower3.5ex\vbox{\hbox{0\rule{0ex}{2.5ex}}
\hbox{\hspace{0.4ex}\rule{.1ex}{1ex}\rule{0ex}{1.4ex}}\hbox{1\strut}}--0--0
\end{tabular}}  & 7 & $\begin{array}{l} (\vp_6,1),(\vp_2,2),(\vp_5,3),
               (2\vp_1,3),\\(\vp_4,4),(\vp_1{+}\vp_7,4),(\vp_1{+}\vp_3,5),\\
        (2\vp_7,5),(\vp_3{+}\vp_7,6),(2\vp_3,7)\end{array}$ & + \\ \hline
$\GR{E}{8}$  & {\small\begin{tabular}{@{}c@{}}
0--1--0--0--\lower3.5ex\vbox{\hbox{0\rule{0ex}{2.5ex}}
\hbox{\hspace{0.4ex}\rule{.1ex}{1ex}\rule{0ex}{1.4ex}}\hbox{0\strut}}--0--0
\end{tabular}}  & 4 & \ $(\vp_1,1),(\vp_7,2),(\vp_2,3),
               (\vp_3,4)$ & +
\\ \cline{2-5}
     & {\small\begin{tabular}{@{}c@{}}
0--0--0--0--\lower3.5ex\vbox{\hbox{0\rule{0ex}{2.5ex}}
\hbox{\hspace{0.4ex}\rule{.1ex}{1ex}\rule{0ex}{1.4ex}}\hbox{1\strut}}--0--0
\end{tabular}}  & 8 &  $\begin{array}{l} (\vp_1,1),(\vp_7,2),(\vp_2,3),
               (\vp_3,4) \\ (\vp_8,4),(\vp_6,5),(\vp_4,6),(\vp_5,7)
               \end{array}$ & +
\\ \hline
$\GR{F}{4}$ & $0{-}0{\Leftarrow} 1{-}0$
      & 4 & \ $(\vp_4,1),(2\vp_1,2),(\vp_3,3),
               (2\vp_2,4)$ & $+$ \\ \hline
$\GR{G}{2}$ & $1{\Lleftarrow} 0\phantom{-0}$ & 2 &
          \ $(\vp_1,1),(\vp_2,1)$ & $-$  \\ \hline
\end{tabular}
\end{center}
\vskip-1ex
\end{table}
{\bf Example.} Take the second nilpotent orbit for $\GR{E}{7}$ in Table~2.
Its dimension is 70 and the Dynkin-Bala-Carter label is $4\GR{A}{1}$. Here
the generators of
$\cs(\g\langle 2\rangle\oplus\g\langle 3\rangle)^{U(L)}$ are:
\\[.7ex]
$f_1=(\vp_6,1),\, f_2=(\vp_2,2),\, f_3=(\vp_5+\vp_1-\vp_7,2),\,
f_4=(\vp_4+\vp_1-\vp_7,3),\, f_5=(2\vp_1,3),\, f_6=(\vp_7-\vp_1,1),\,
f_7=(\vp_3-\vp_7,1)$.
Therefore
$\cs(\g\langle 2\rangle\oplus\g\langle 3\rangle)^{U(L)}_{(+)}$ is
generated by \\[.7ex]
\centerline{$f_1,\, f_2,\, f_3f_6,\, f_5,\, f_4f_6,\, f_5f_6,\, f_5f_6f_7,\,
f_5f_6^2,\, f_5f_6^2f_7,\, f_5f_6^2f_7^2$\ .}
\\[.6ex]
This ordering of generators corresponds
to their ordering in Table~2.
\\[.7ex]
Comparing our computations with known structure of $\bbk[G{\cdot}e]$
as $G$-module when possible, we see that Conjecture~\ref{conj1prim}
is valid for several items in Table~2; that is, the generators of
$\cs(\g\langle 2\rangle\oplus\g\langle 3\rangle)^{U(L)}_{(+)}$ given in the
fourth column do belong to $\bbk[G{\cdot}e]^U$.
As far as I know, this refers to the orbits in
$\GR{E}{6},\GR{F}{4},\GR{G}{2}$, and the second orbit for $\GR{E}{8}$.
The last case is especially interesting. It was proved in \cite{e8} that
the algebra of regular functions on this 128-dimensional orbit is a model
algebra for $\GR{E}{8}$, i.e., each simple finite-dimensional
$\GR{E}{8}$-module occurs in $\bbk[G{\cdot}e]$ exactly once.
On the other hand,
$\cs(\g\langle 2\rangle\oplus\g\langle 3\rangle)^{U(L)}_{(+)}$ is a polynomial
algebra and the weights of free generators are exactly the fundamental
weights of $\GR{E}{8}$. As $\hat\tau^o$ is an embedding of
$\bbk[G{\cdot}e]^U$ into
$\cs(\g\langle 2\rangle\oplus\g\langle 3\rangle)^{U(L)}_{(+)}$,
we see that
$\hat\tau^o$ is an isomorphism in this case.
Furthermore, our computation provides the degree in which each fundamental
$\GR{E}{8}$-module appears in the graded algebra $\bbk[\ov{G{\cdot}e}]^U$.
\\[.6ex]
It is instructive to observe that
$\cs(\g\langle 2\rangle\oplus\g\langle 3\rangle)^{U(L)}$
is always polynomial, whereas this is not always the case for the
subalgebra $\cs(\g\langle 2\rangle\oplus\g\langle 3\rangle)^{U(L)}_{(+)}$.
For instance, the latter has $(t+1)(t+4)/2$ generators in the case of
${\frak so}_{4t+4}$ or ${\frak so}_{4t+l+3}$,
whereas its
Krull dimension equals $2t+2$. Modulo Conjecture~\ref{conj1prim}, this means
that the homological dimension of the algebra of covariants for the
closure of a spherical nilpotent orbit may be arbitrarily large.
\\[.7ex]
{\bf Remark.} Using condition $(\spadesuit)$ given in the proof of
Theorem~\ref{height2}(iii), our computations for the orbits of height 3
yield the following assertions:
\begin{itemize}
\item If $\hot(e)=3$, then
$\cs(\g\langle 2\rangle\oplus\g\langle 3\rangle)$ is a free
$\cs(\g\langle 2\rangle\oplus\g\langle 3\rangle)^{U(L)}$-module.
\item the validity of Conjecture~\ref{conj1prim} would imply that if
$\bbk[G{\cdot}e]^U$ is polynomial, then $\bbk[G{\cdot}e]$ is
a free $\bbk[G{\cdot}e]^U$-module.
\end{itemize}

\end{document}